\documentclass[%
11pt,
onecolumn,
tightenlines,
superscriptaddress,
preprintnumbers,
nofootinbib,
amsmath,amssymb,amsthm,
physrev,
eqsecnum,
]{revtex4-2}

\usepackage{isomath}
\usepackage{amsmath,amsthm}
\usepackage{amsbsy}
\usepackage{amssymb}
\usepackage{amscd}
\usepackage{amsfonts}
\usepackage{stmaryrd}
\usepackage{siunitx}
\usepackage{euscript}
\usepackage[utf8]{inputenc}
\usepackage[T1]{fontenc}
\usepackage{newtxtext} 
\everymath{\displaystyle}
\usepackage{exscale}

\usepackage{graphicx}
\usepackage{boxedminipage}
\usepackage{calc}
\usepackage[usenames,dvipsnames]{xcolor}
\graphicspath{ {media/} }
\usepackage[caption=false,justification=centerlast]{subfig}

\usepackage{setspace}
\usepackage{enumitem}
\setitemize{noitemsep,topsep=0pt,parsep=0pt,partopsep=0pt}
\setenumerate{noitemsep,topsep=0pt,parsep=0pt,partopsep=0pt}
\setdescription{noitemsep,topsep=0pt,parsep=0pt,partopsep=0pt}

\usepackage{hyperref}
\usepackage[normalem]{ulem}

\usepackage[small]{titlesec}

\titlespacing*{\section}{0pt}{12pt plus 4pt minus 2pt}{2pt plus 2pt minus 2pt}
\titlespacing*{\subsection}{0pt}{12pt plus 4pt minus 2pt}{2pt plus 2pt minus 2pt}
\titlespacing*\subsubsection{0pt}{12pt plus 4pt minus 2pt}{2pt plus 2pt minus 2pt}
\titlespacing*\paragraph{0pt}{12pt plus 4pt minus 2pt}{2pt plus 2pt minus 2pt}

\makeatletter
    \renewcommand*{\thesection}{\arabic{section}}
    \renewcommand*{\thesubsection}{\thesection.\Alph{subsection}}
    \renewcommand*{\p@subsection}{}
    \renewcommand*{\thesubsubsection}{\thesubsection.\arabic{subsubsection}}
    \renewcommand*{\p@subsubsection}{}
\makeatother

\usepackage{isomath}
\usepackage{amsmath}
\usepackage{amssymb}
\usepackage{amscd}
\usepackage{amsfonts}

\newcommand{\calL}{{\cal L}}

\newcommand{\calS}{{\cal S}}

\newtheorem{remark}{Remark}[section]

\newcommand{\abs}[1]{\lvert#1\rvert}

\newcommand{\parderiv}[2]{\frac{\partial #1}{\partial #2}}
\newcommand{\dm}{\ \mathrm{d}}
\newcommand{\deriv}[2]{\frac{\dm #1}{\dm #2}}

\newcommand{\bfk}{{\mathbold k}}

\newcommand{\bfn}{{\mathbold n}}

\newcommand{\bfx}{{\mathbold x}}

\newcommand{\bfA}{{\mathbold A}}

\newcommand{\bfD}{{\mathbold D}}

\newcommand{\bfK}{{\mathbold K}}

\newcommand{\bfN}{{\mathbold N}}

\newcommand{\bfP}{{\mathbold P}}
\newcommand{\bfQ}{{\mathbold Q}}

\usepackage{siunitx}
\usepackage{soul}
\DeclareMathOperator{\RE}{Re}
\DeclareMathOperator{\IM}{Im}

\begin{document}


\preprint{Accepted to appear in Journal of the Mechanics and Physics of Solids (DOI: \url{https://doi.org/10.1016/j.jmps.2022.104992})}

\title{\Large{
    Multiband Homogenization of Metamaterials in Real-Space:
    \\
    Higher-Order Nonlocal Models and Scattering at External Surfaces 
}}

\author{Kshiteej Deshmukh}
    \email{kjdeshmu@andrew.cmu.edu}
    \affiliation{Department of Civil and Environmental Engineering, Carnegie Mellon University}

\author{Timothy Breitzman}
    \affiliation{Air Force Research Laboratory}
    
\author{Kaushik Dayal}
    \affiliation{Department of Civil and Environmental Engineering, Carnegie Mellon University}
    \affiliation{Center for Nonlinear Analysis, Department of Mathematical Sciences, Carnegie Mellon University}
    \affiliation{Department of Mechanical Engineering, Carnegie Mellon University}
    
\date{\today}


\begin{abstract}
	This work develops a dynamic linear homogenization approach in the context of periodic metamaterials.
	By using approximations of the dispersion relation that are amenable to inversion to real-space and real-time, it finds an approximate macroscopic homogenized equation with constant coefficients posed in space and time; however, the resulting homogenized equation is higher order in space and time.
	The homogenized equation can be used to solve initial-boundary-value problems posed on arbitrary non-periodic macroscale geometries with macroscopic heterogeneity, such as bodies composed of several different metamaterials or with external boundaries.
    The approach is applied here to problems with scalar unknown fields in one and two spatial dimensions.
    
	First, considering a single band, the dispersion relation is approximated in terms of rational functions, enabling the inversion to real space.
	The homogenized equation contains strain gradients as well as spatial derivatives of the inertial term.
	Considering a boundary between a metamaterial and a homogeneous material, the higher-order space derivatives lead to additional continuity conditions. 
	The higher-order homogenized equation and the continuity conditions provide predictions of wave scattering in 1-d and 2-d that match well with the exact fine-scale solution; compared to alternative approaches, they provide a single equation that is valid over a broad range of frequencies, are easy to apply, and are much faster to compute.

	Next, the setting of two bands with a bandgap is considered.
	The homogenized equation has also higher-order time derivatives.
    Notably, the homogenized model provides a single equation that is valid over both bands and the bandgap.
    The continuity conditions for the higher-order spatio-temporal homogenized equation are applied to wave scattering at a boundary, and show good agreement with the exact fine-scale solution.
    The method is also applied to a problem with multiple scattered propagating waves for which the classical jump conditions cannot provide even approximate solutions, and the results are shown to match reasonably well with the exact fine-scale solutions.
    
	Using that the order of the highest time derivative is proportional to the number of bands considered, a nonlocal-in-time structure is conjectured for the homogenized equation in the limit of infinite bands.
	This suggests that homogenizing over finer length and time scales -- with the temporal homogenization being carried out through the consideration of higher bands in the dispersion relation -- is a mechanism for the emergence of macroscopic spatial and temporal nonlocality, with the extent of temporal nonlocality being related to the number of bands considered.
    
\end{abstract}

\maketitle


\section{Introduction}

\begin{figure}[htb!]
    \centering
    \includegraphics[width=0.8\textwidth]{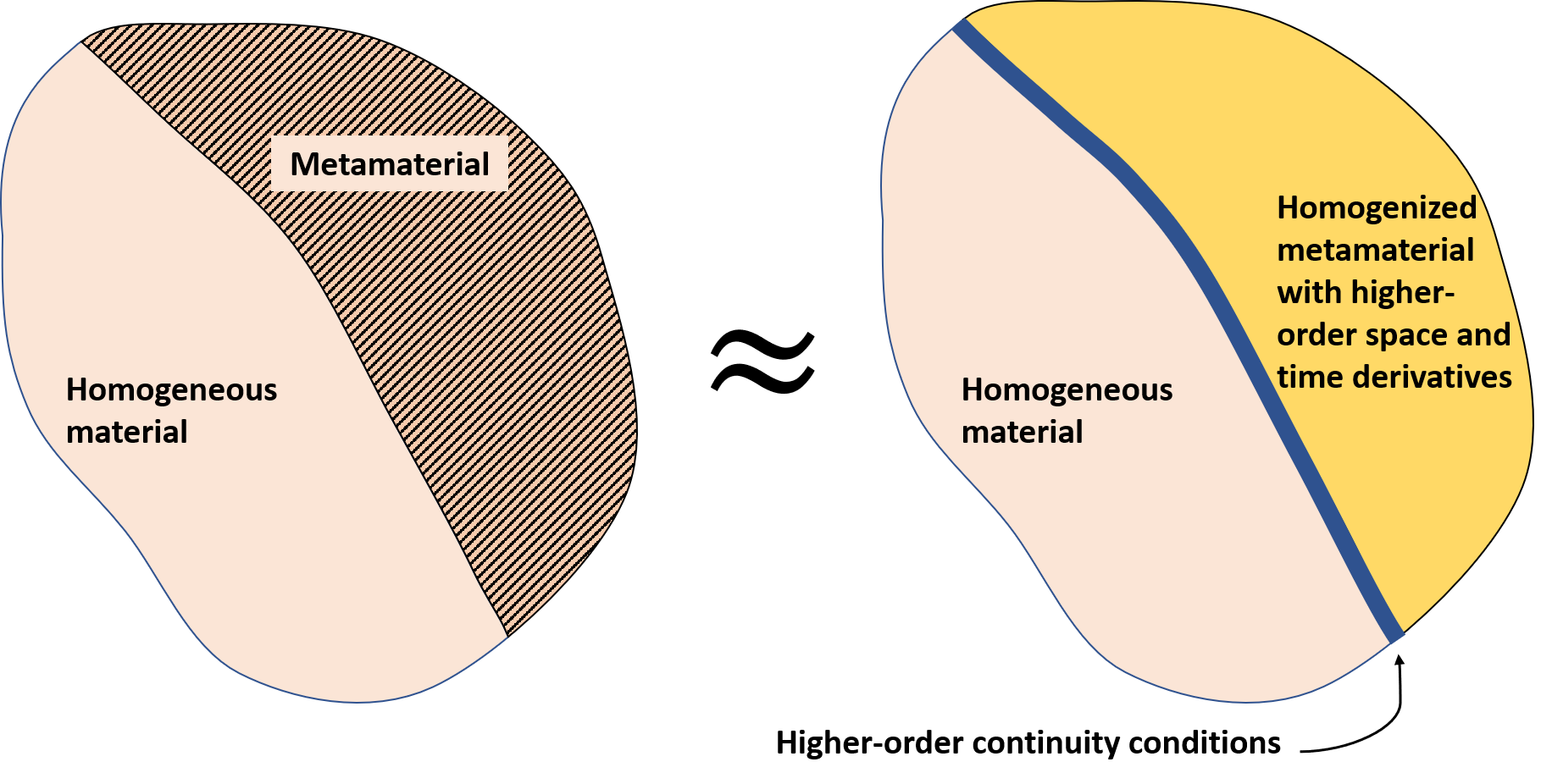}
    \caption{The metamaterial (left) gives a homogenized material with higher-order space-time derivatives (right), and consequent higher-order continuity conditions.
    }
    \label{fig:method_sketch}
\end{figure}

Metamaterials have the potential to display unusual properties, and have been the focus of much attention in mechanics and electromagnetism, e.g. \cite{Deng2021,liu2018broadband,shmuel2012band,Hussein2007,Matlack2018,Dayal2007a,khajehtourian2021continuum,Bertoldi2017,gazonas2006genetic,spadoni2009phononic,hussein2014dynamics,lipton2018effective,filipov2015origami,hajarolasvadi2021dispersion,eugster2019continuum,giorgio2019two} and numerous others.
A primary focus of much of this effort is to find the dispersion relations in the context of infinite periodic systems, and then tailor or optimize these relations to obtain desirable properties.

In this paper, we use the dispersion relation to develop homogenized models, in the context of linear periodic metamaterials.
The homogenized models that we develop have two key features:
(1) they are posed in real space and time, i.e., not in frequency space; and
(2) they are posed in terms of a single equation that is valid over a large frequency range and across multiple bands.
These features can enable the application of the models to problems posed on complex geometries and with complex time-dependent behavior.

Our approach is based on developing approximations to the dispersion relation and then inverting it to obtain the real-space model.
This requires us to balance between two competing objectives: (1) an accurate representation of the dispersion relation, and (2) tractability of the representation to inversion to real space to obtain standard differential operators.
For a single band, we find that a class of models that go back to Toupin \cite{toupin1962elastic} provides a good balance.
In brief, this class of models can be considered as approximating dispersion curves as rational functions; the inversion to real-space then provides standard but higher-order differential operators.
Specifically, the homogenized equation includes strain gradients as well as spatial gradients of the inertial terms.
When we extend this to higher bands, we find the appearance of higher-order time derivatives as well.

To summarize our overall procedure, we: (1) start with a given dispersion relation that is derived under assumptions of periodicity; (2) use a rational function approximation of the dispersion relation to develop a homogenized equation in space and time with constant coefficients; and (3) then {\em assume} that the homogenized equation is still valid even when the coefficients vary in space and the body is finite.
Roughly, our procedure assumes implicitly a separation of scales between the periodic microstructure and the macroscopic heterogeneity.

To test and demonstrate the homogenized model, we turn to an important topic of current interest in the community: to predict the scattering of waves at boundaries between different metamaterials or boundaries between a metamaterial and a homogeneous material, e.g. \cite{Srivastava2017,willis2020transmission-JMPS,willis2020transmission-PRS,Cakoni2019}.
We use problems extracted from this body of work as examples on which to test our approach, and discuss this in detail below.
Our approach, in brief, is to identify the additional nonstandard continuity conditions that arise at a surface of discontinuity due to the higher-order derivatives in the homogenized equations.
These higher-order continuity conditions then provide precisely the information required to uniquely solve for the wave scattering at the boundary.
A schematic view of this approach is shown in Figure \ref{fig:method_sketch}.

\paragraph*{Prior Work on Scattering at Boundaries between Metamaterials and Homogeneous Materials.}

A simplistic approach to model the scattering at a boundary might start from the usual assumption of infinite periodicity to obtain the homogenized properties of the materials on either side, and then use these properties in the nonperiodic setting with a boundary.
However, as pointed out in \cite{Srivastava2017}, such an approach would consider only the propagating waves and miss the contribution from the evanescent waves that are present in the nonperiodic setting.
They showed that if such homogenized models are applied to this problem, the evanescent modes cannot be accounted for, and the classical continuity conditions -- continuity of displacement and traction -- must be satisfied with propagating waves only; this leads to a violation of the energy flux conservation. 

\cite{Srivastava2017} discuss an approach to deal with this problem.
Their approach is to minimize the errors in displacement continuity and traction continuity, subject to the conservation of energy, by considering only the propagating modes.
While simple and easy to use, it does not consider the microstructural information in formulating this criterion, i.e., it does not account for the details of the evanescent modes.
Further, it can be difficult to apply this in a more general setting, for instance in a setting with complex geometry and time-dependence where a decomposition into propagating modes is difficult.

A different ``bottom-up'' approach to this problem, and analogous problems, has been proposed by other workers, e.g., \cite{Maurel2018a, Cornaggia2020, Cakoni2019, pham2021revisiting,MarigoMaurel2017,Maurel2018,Maurel2019,Guzina2019,Meng2018, marigo2016two}.
These approaches rigorously account for the lack of periodicity, e.g., by introducing boundary-corrector functions or using matched asymptotic expansions in 2-scale methods. 
Such methods typically require the coupled solution of the cell problem -- solving for the fast variable in a periodic unit cell -- along with the effective PDE describing the wave transmission, or the computation of the boundary-corrector functions.
These requirements make the approaches challenging and computationally expensive, and consequently difficult to apply in a general setting.
Further, these methods typically are restricted to a relatively limited frequency band, making them difficult to apply to general problems with complex time-dependent loading.

\paragraph*{The Proposed Approach.}

The approach that we propose aims to sit between the approach of \cite{Srivastava2017} and the rigorous homogenization approaches.
Specifically, we aim for our homogenized models to be (1) applicable to a broad range of frequencies, (2) applicable to general geometries and time-dependent loadings, (3) computationally efficient and easy to compute with.

In the specific context of the insights provided by \cite{Srivastava2017} on the importance of evanescent waves, our homogenized models are able to capture the evanescent modes and the bandgap.
Specifically, in higher dimensions, there are typically 2 propagating scattered waves: a reflected wave and a transmitted wave.
While the classical continuity conditions -- i.e., that traction and displacement are continuous --  enable us to approximately find these propagating waves, this approach neglects completely the effects of the evanescent modes. 
Further, \cite{Srivastava2017} show that there exist important regimes wherein there are more than 2 propagating waves, and the 2 classical continuity conditions are insufficient to approximately find the propagating waves.

Our approach provides higher-order continuity conditions (Fig. \ref{fig:method_sketch}) at the interface that augment the standard continuity conditions.
In general, these additional continuity conditions enable us to better account for the evanescent waves, thereby providing more accurate solutions.
In the specific situations studied by \cite{Srivastava2017}, the additional conditions enable us to uniquely find the propagating scattered waves.

\paragraph*{Emergence of Temporal Nonlocality with Multiple Bands.} 

Existing dynamic homogenization techniques can provide a homogenized equation for a specific band of the dispersion curve by expanding about the frequency in that band, e.g. \cite{Guzina2019, Craster2010,Harutyunyan2016}. 
However, this typically does not provide a single equation that covers several bands.
In this work, we propose a single homogenized equation that is applicable over a range of frequencies that covers multiple bands and the bandgaps.
An important feature of this homogenized equation is that it has higher derivatives in time, suggesting the emergence of nonlocality in time due to homogenization as we increase the number of bands.
Other researchers have developed the idea of nonlocality in time as an outcome of homogenization: foundational early work in this area includes the approach by Willis \cite{Willis1980, willis1981variational, willis1997dynamics}, and the mathematical work of Tartar \cite{tartar1989nonlocal,tartar1991memory} and others following him \cite{antonic1993memory}.
However, these and other prior works does not seem to have explored the path of accounting for higher-order dispersion bands as a mechanism for temporal nonlocality.

\paragraph*{Organization.}

In Section \ref{sec:fine_cale}, we formulate the exact fine-scale model and discuss the specific examples that we will use to compare the exact and homogenized models.
In Sections \ref{sec:1d} and \ref{sec:2d}, we describe the development of the homogenized single-band model in 1-d and 2-d respectively, and compare its predictions with those of the exact model.
In Section \ref{sec:OA}, we describe the development of the homogenized model for two bands with a bandgap, and compare its predictions with those of the exact model.
In Section \ref{sec:2D-2Band}, we apply our approach to the setting in which there can be more than 2 propagating waves, wherein the classical continuity conditions are insufficient to provide a unique solution.
In Section \ref{sec:integral_models}, we discuss nonlocal equations as possible homogenized models that account for an infinite number of bands.

\section{Formulation and Fine-Scale Model}
\label{sec:fine_cale}
\subsection{Notation}

We will refer to the model that considers the detailed microstructure as the fine-scale model, and its solution as the fine-scale solution.
Analogously, we will refer to the homogenized effective dynamical PDE as the homogenized model, and its solution as the homogenized solution.

The primary variable, i.e. the displacement in the mechanical setting, will be denoted $u(\bfx,t)$, where $t$ denotes time and $\bfx$ is the spatial coordinate.
In 1-d, $\bfx \equiv x$; in 2-d, $\bfx \equiv (x_1, x_2)$ in Cartesian coordinates.
Our 2-d problems will be in the antiplane setting, so $u$ is always a scalar.

For conciseness, we use the following representation for derivatives: 
$\partial_x \equiv \parderiv{}{x}$, $\partial_t \equiv \parderiv{}{t}$, $\partial_x^2 \equiv \frac{\partial^2}{\partial x^2}$, $\partial_t^2 \equiv \frac{\partial^2}{\partial t^2}, \partial_i \equiv \parderiv{}{x_i}$ and so on.
We use the summation convention, e.g., $\partial_{ii} \equiv \partial^2_1 + \partial^2_2$.

We consider bodies wherein the material properties -- either in the fine-scale model or in the homogenized model -- can be discontinuous across boundaries.
Consequently, various quantities can be discontinuous at these boundaries, and the jump in a quantity across a surface of discontinuity is denoted $\llbracket \cdot \rrbracket$.
The surfaces of discontinuity are assumed to be fixed in space.

The Bloch wavevector in the fine-scale models is denoted $\bfK = (K_1, K_2)$, and the wavevector in the homogenized setting is denoted $\bfk = (k_1, k_2)$. The frequency is denoted $\omega$, and will be normalized by the midgap frequency $\omega_0$.

\subsection{Fine-scale Governing Equations} 

We will consider the linear wave equation in heterogeneous media.
Denoting the domain of the body by $\Omega$, we have:
\begin{align}
    \text{1-d: } \quad \rho(x) \partial_{t}^2 u(x,t) &= \partial_{x} \left( E(x) \partial_x u(x,t) \right) \quad \text{ in } \Omega
    \\
    \text{2-d antiplane: } \quad \rho(\bfx) \partial_{t}^2 u(\bfx,t) &= \partial_{i} \left( \mu(\bfx) \partial_{i} u(\bfx,t)\right) \quad \text{ in } \Omega
\end{align}
where the fine-scale material is assumed to be isotropic for antiplane shear, and $\mu(\bfx)$ and $E(x)$ denote moduli, and $\rho$ is the density.

Let $\mathcal{S}$ denote a stationary surface of discontinuity in $\Omega$.
The continuity conditions on $\mathcal{S}$ are given by:
\begin{align}
      \text{1-d: } \quad & \left\llbracket u(x,t) \right\rrbracket =0, \quad \left\llbracket E(x) \partial_x u(x,t) \right\rrbracket =0 \quad \forall x \in \mathcal{S}
     \\
    \text{2-d antiplane:}  \quad & \left\llbracket u(\bfx,t) \right\rrbracket =0, \quad \left\llbracket \mu(\bfx) \partial_{i} u(\bfx,t) \right\rrbracket \hat n_i(\bfx) =0  \quad \forall \bfx \in \mathcal{S}
\end{align}
where $\hat\bfn(\bfx)$ is the unit normal to $\mathcal{S}$. 
The continuity conditions correspond to the classical conditions of displacement continuity and traction continuity respectively.

\subsection{Problem Geometries and Boundary Conditions}
\label{sec:formulation-problems}

\paragraph{One-dimensional Problems. }

\begin{figure}[ht!]
	\subfloat[
	    Single boundary geometry in 1-d, with a semi-infinite laminate metamaterial adjoining a  semi-infinite homogeneous material.
	    \label{fig:single_interface_setup}
    ]
	{\includegraphics[width=0.44\textwidth]{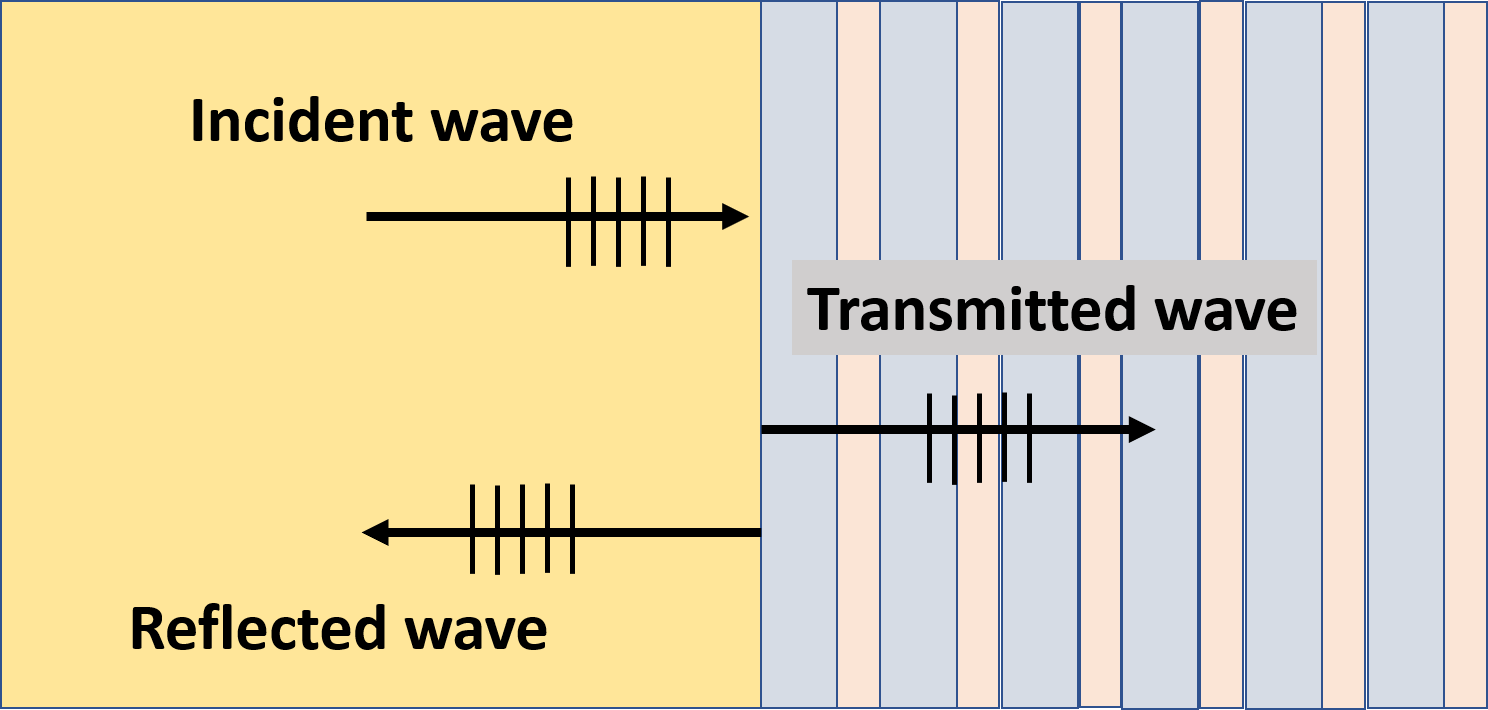}}
	\hfill
    \subfloat[
        Double boundary geometry in 1-d, with a laminate metamaterial of finite extent $L$ (taken to be $10$ unit cell widths) between two identical semi-infinite homogeneous materials.
        \label{fig:2interface_setup}
    ]
    {\includegraphics[width=0.46\textwidth]{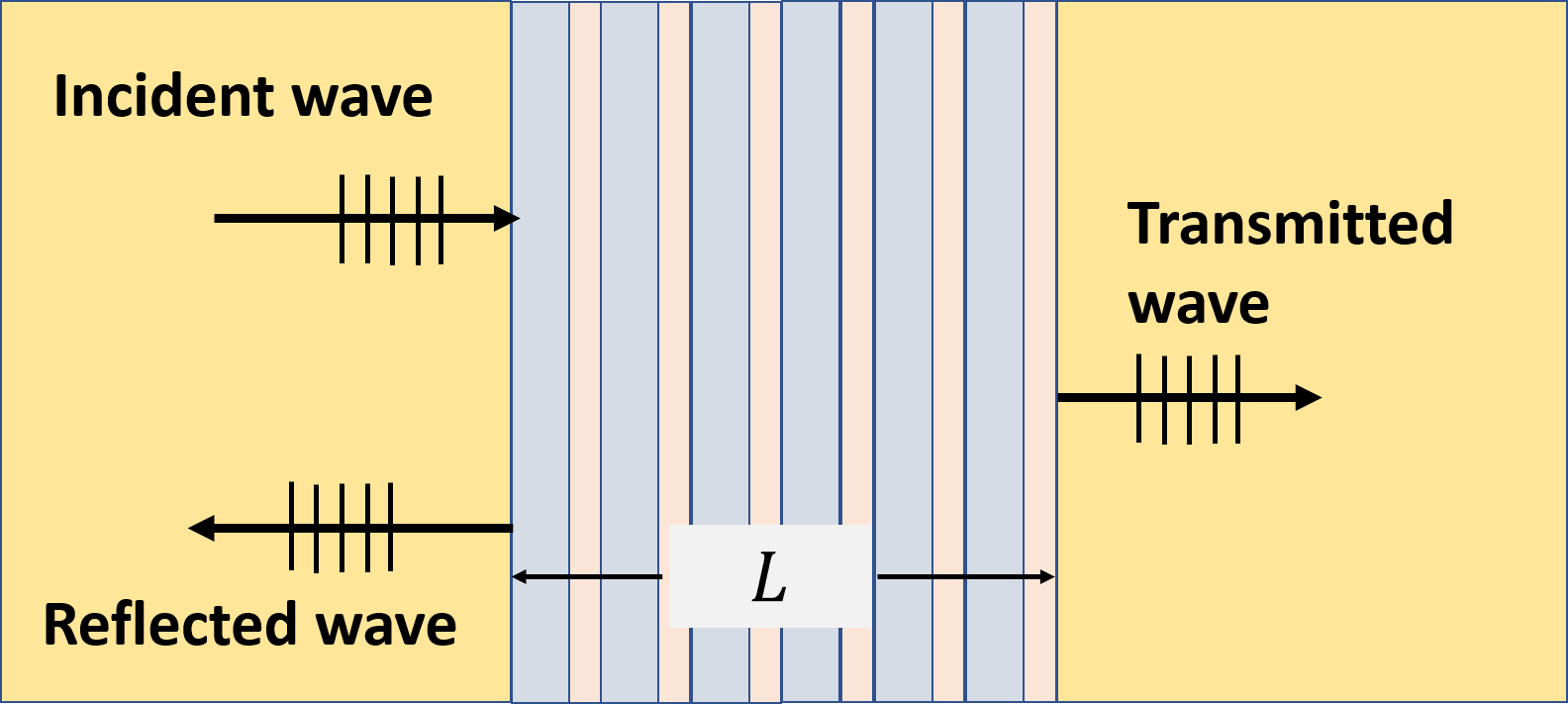}}
   \caption{Schematics of the problem geometries in 1-d.}
   \label{fig:1d_setups}
\end{figure}

In 1-d, we consider two problems, as shown in Figure \ref{fig:1d_setups}.
In both problems, a given traveling wave originates at $x\to-\infty$ in the homogeneous material and is incident at the boundary between the homogeneous material and the metamaterial.
At $x\to\infty$, a Bloch wave form is used in the semi-infinite metamaterial, and a traveling wave is used in the semi-infinite homogeneous material in the single- and double- boundary problems respectively.
This formulation, and the solution, can be found in many works in the literature, e.g., Section $2$ in \cite{Willis2016} and \cite{Cornaggia2020}.

The primary quantities of interest are the magnitudes and energies of the reflected and transmitted waves at steady state; we aim to compare the homogenized model predictions of these quantities with the fine-scale model predictions for various frequencies.

\paragraph{Two-dimensional Problem. }

\begin{figure}[ht!]
    \centering
    \includegraphics[width=0.5\textwidth]{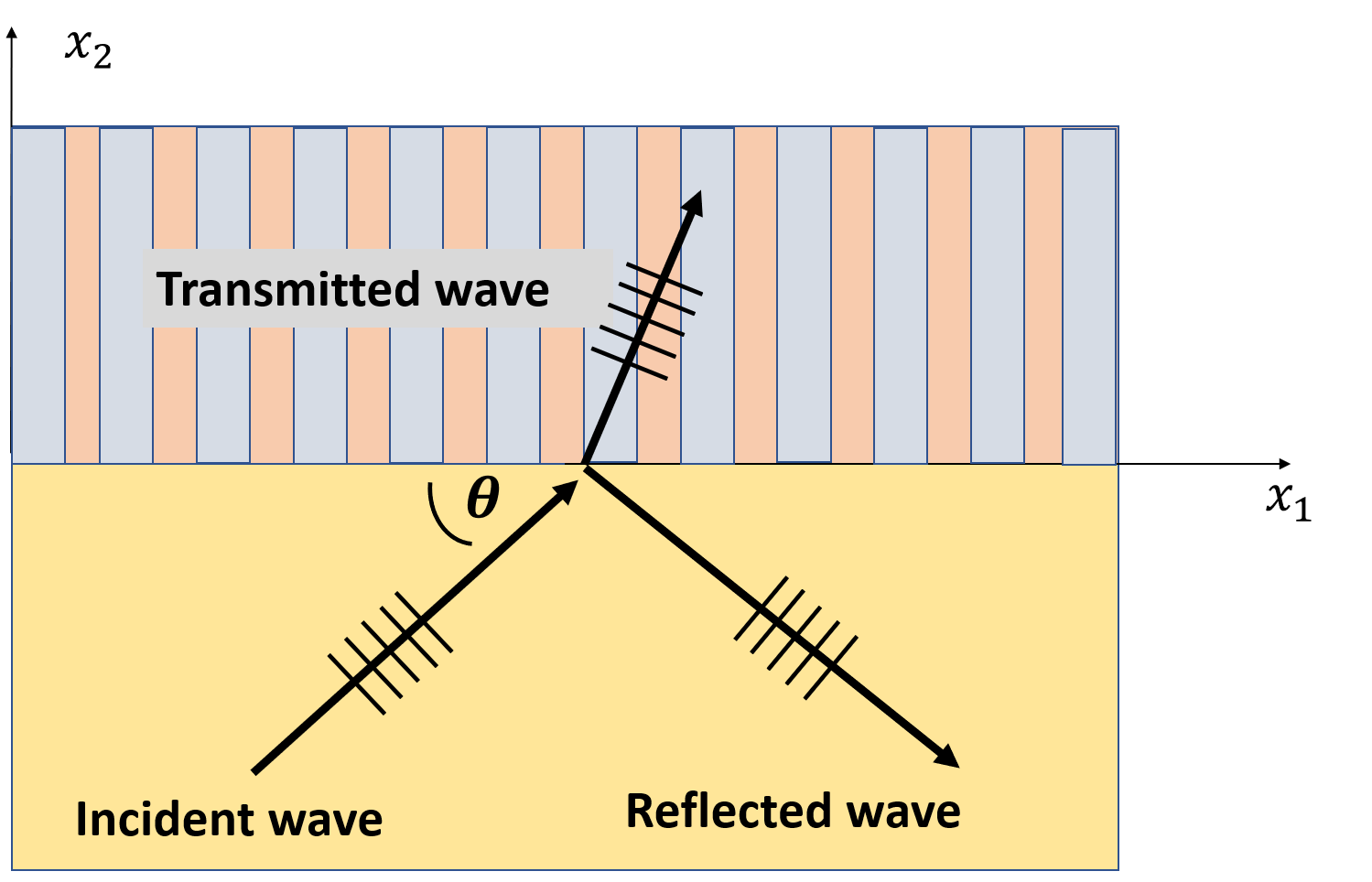}
    \caption{Problem geometry in 2-d, with a semi-infinite laminate metamaterial adjoining a semi-infinite homogeneous material. The macroscopic boundary is normal to the laminae in the metamaterial.}
    \label{fig:perp_interface_setup}
\end{figure}

In 2-d, we consider the geometry shown shown in Figure \ref{fig:perp_interface_setup}.
We consider plane traveling waves originating at infinity in the semi-infinite homogeneous material, and oriented at an angle $\theta$ to the tangent of boundary of the metamaterial.
The plane wave scatters on the boundary giving rise to transmitted and reflected waves.
Because the boundary is infinite and planar, the transmitted and reflected waves are also plane waves and we apply a Bloch form at infinity in the semi-infinite metamaterial.
This formulation, and the solution, follows \cite{Willis2016, Srivastava2017}.

The primary quantities of interest are the magnitudes and energies of the transmitted and reflected waves at steady state, and we compare the homogenized model predictions and fine-scale model predictions for various frequencies and various angles of incidence $\theta$.

\subsection{Unit Cell Material Properties}
\label{sec:num-properties}

\begin{figure}
    \centering
    \includegraphics[width=0.3\textwidth]{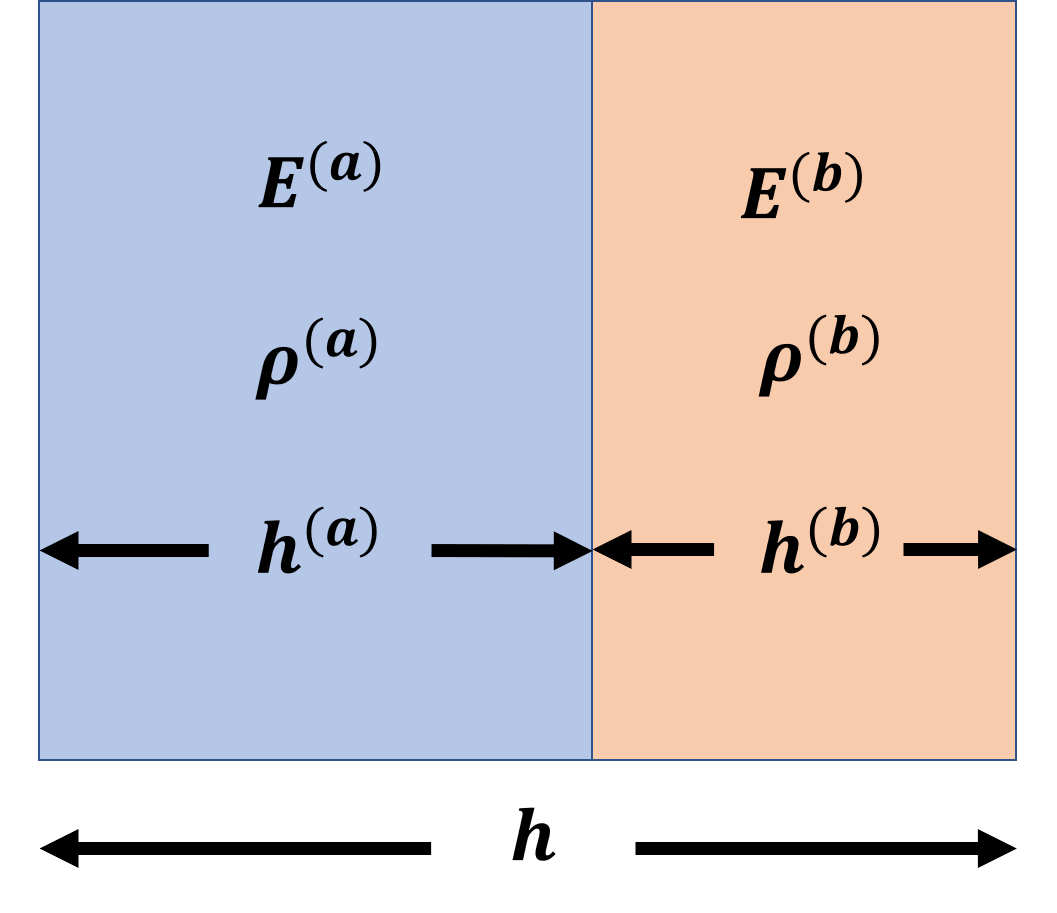}
    \caption{Unit cell for laminate metamaterial. $E$ denotes modulus and $\rho$ denotes density. }
    \label{fig:unit_cell}
\end{figure}

Throughout this paper, we consider a periodic bilayer laminate metamaterial composed of two materials $a$ and $b$; the geometry and notation is shown in Figure \ref{fig:unit_cell}.
Table \ref{tab:composite_prop1d} lists the numerical values of the material properties chosen for all the examples (both 1-d and 2-d) involving a single band.
Table \ref{tab:composite_prop1d_bettergap} lists the numerical values of the material properties chosen for the examples  involving 2 bands in Section \ref{sec:OA}; the key difference is a bigger contrast to have a larger bandgap.  
The bandgap frequencies used for normalization for the properties from  Table \ref{tab:composite_prop1d} and Table \ref{tab:composite_prop1d_bettergap} are  $\omega_0=\SI{3.448685}{\radian\per\second}$ and $\omega_0 =\SI{8.64675}{\radian\per\second} $ respectively.

The homogeneous material has modulus and density denoted $E_h$ and $\rho_h$ respectively. 
We use a numerical value of $\rho_h =1$ throughout.
For the 1-d single boundary problem, we use $E_h=0.25$ when we examine 1-band approximations and $E_h=8$ when we examine 2-band approximations; for the 1-d double boundary problem and the 2-d problem, we use $E_h=2$.
These values are chosen to provide numerical results that can be easily compared when plotted.

\begin{table}[ht!]
    \parbox{.45\textwidth}{
    \begin{tabular}{|l|l|}
       \hline
       Layer &  $a/b$ \\
       \hline 
       modulus: $E^{(a)}/E^{(b)}$ & $1/5$ \\
       \hline 
       thickness: $h^{(a)}/h^{(b)}$   & $4$ \\
       \hline
       density:$\rho^{(a)}/\rho^{(b)}$ & $1$ \\
        \hline
    \end{tabular}
    \caption{Properties of the periodic laminate considered for all calculations (1-d and 2-d) involving a single band, with \{$E^{a}$, $\rho^{(a)}$, $h^{(a)}$\} = \{$1,1,4/5$\}.}
    \label{tab:composite_prop1d}
    }
    \hspace{0.05\textwidth}
    \parbox{.45\textwidth}{
    \begin{tabular}{|l|l|}
        \hline
       Layer &  $a/b$ \\
       \hline 
       modulus: $E^{(a)}/E^{(b)}$ & $1/20$ \\
       \hline 
       thickness: $h^{(a)}/h^{(b)}$   & $1/6$ \\
       \hline
       density:$\rho^{(a)}/\rho^{(b)}$ & $1$ \\
        \hline
    \end{tabular}
    \caption{Properties of the periodic laminate for the multiband calculations in Section \ref{sec:OA}, with \{$E^{a}$, $\rho^{(a)}$, $h^{(a)}$\} = \{$1,1,1/7$\}.}
    \label{tab:composite_prop1d_bettergap}
    }
\end{table}

\subsection{Dispersion Relations for the Fine-scale Model} \label{sec:composite_charac}

\begin{figure}[ht!]
	\subfloat[
	    Dispersion relation for the laminate in 1-d with the properties from Table \ref{tab:composite_prop1d}, restricted to the first Brillouin zone. 
	    \label{fig:f_vs_K1_1d}
    ]
	{\includegraphics[width=0.47\textwidth]{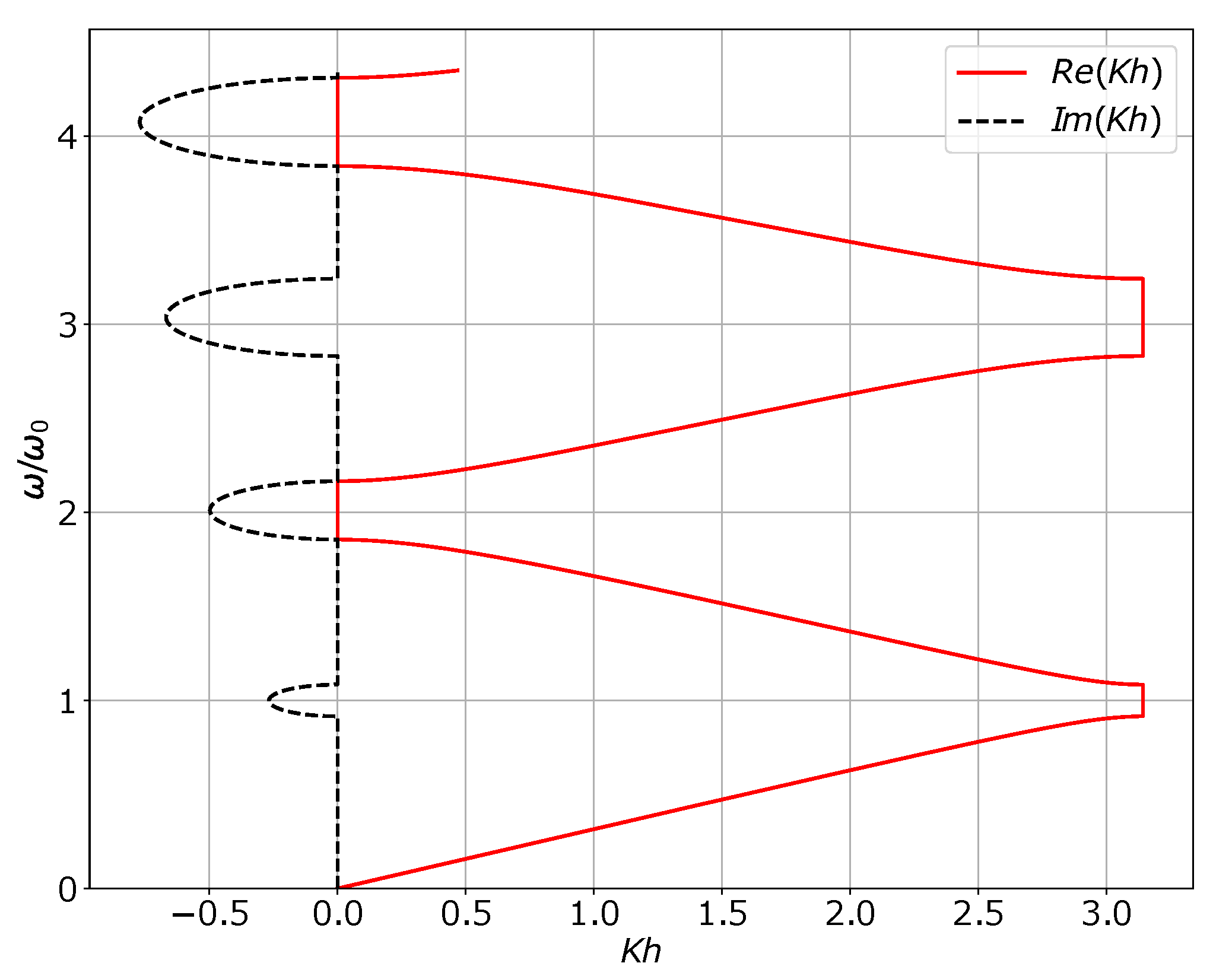}}
    \hfill    
    \subfloat[
        Dispersion relation for the laminate in 2-d with the properties from Table \ref{tab:composite_prop1d}.
        \label{fig:3d_RF22_approx}
    ]
    {\includegraphics[width=0.5\textwidth]{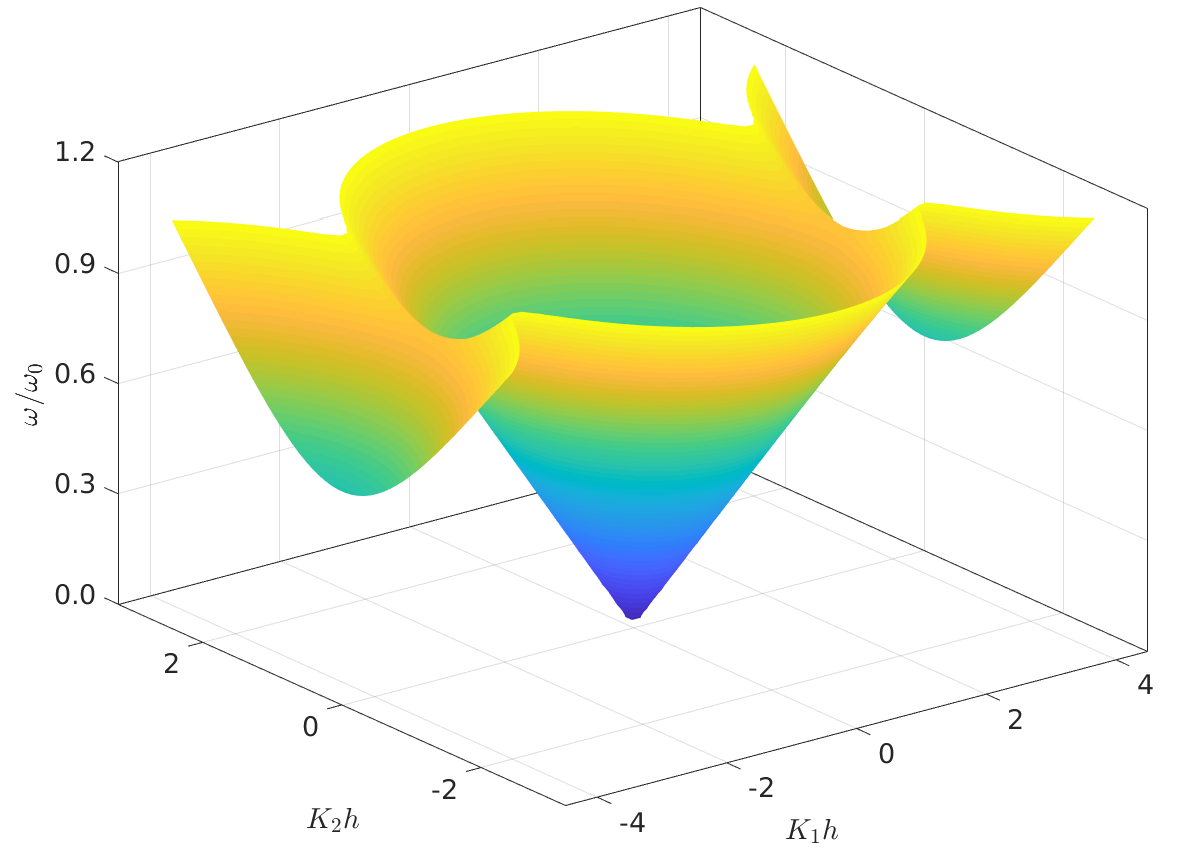}}
    \caption{Dispersion relations in 1-d and 2-d.}
    \label{fig:dispersion-relations}
\end{figure}

The dispersion relation for a periodic laminate in 1-d is classical and given in, e.g., \cite{Willis2016, lekner1994light}:
\begin{equation}\label{disp_rel_periodic_laminate_1d}
    K = \frac{1}{h} \arccos{\left[ \cos{k^{(b)}h^{(b)}}\cos{k^{(a)}h^{(a)}} - \frac{1}{2} \left(\frac{E^{(a)}k^{(a)}}{E^{(b)}k^{(b)}} + \frac{E^{(b)}k^{(b)}}{E^{(a)}k^{(a)}}\right) \sin{k^{(a)}h^{(a)}} \sin{k^{(b)}h^{(b)}} \right]}
\end{equation}
Here, $K$ is the Bloch wavevector, $k^{(a)} = \omega / \sqrt{E^{(a)}/\rho^{(a)}}, k^{(b)} = \omega / \sqrt{E^{(b)}/\rho^{(b)}}, h= h^{(a)} + h^{(b)} $. 
Figure \ref{fig:f_vs_K1_1d} plots this dispersion relation using the properties in Table \ref{tab:composite_prop1d}. 

The dispersion relation $\omega(\bfK)$ for a periodic laminate in 2-d is given in, e.g., \cite{Willis2016, Cornaggia2020}.
It is plotted in Figure \ref{fig:3d_RF22_approx} as a function of the wavevector components $\bfK = (K_1, K_2)$.

\section{Single Band Homogenized Model in One Dimension}
\label{sec:1d}

In this section, we examine the simple setting of a 1-d model with a single band.
We begin by identifying an appropriate approximation to the dispersion relation; the approximation captures the primary feature that waves above a certain frequency are not supported, and at the same time is amenable to be inverted to find the homogenized real-space-time dynamical equation (Section \ref{sec:1d-1band-dispersion}).
We then identify the continuity conditions that correspond to the homogenized dynamical equation in a macroscopically-heterogeneous medium (Section \ref{sec:1d-1band-jump}).
Finally, we apply these continuity conditions to study scattering at macroscopic boundaries in specific examples (Section \ref{sec:1d-1band-examples}). 

\subsection{Approximating the Dispersion Relation}
\label{sec:1d-1band-dispersion}

We begin by comparing the dispersion relations of classical linear elasticity, strain gradient elasticity, and ``microinertia'' models in the context of homogeneous materials.
Such models were introduced by Toupin \cite{toupin1962elastic}, and have developed and applied by numerous other researchers; a sample of this vast literature includes \cite{mindlin1963microstructure,lim2015higher,lam2003experiments,amanatidou2002mixed,zhang2005inclusions,maranganti2007novel,fried2009gradient}.
The terminology of ``microinertia'' is more recent \cite{askes2006gradient,askes2008four,askes2011gradient,hui2014high}, and refers to terms that have spatial derivatives of second-order time derivatives.
\begin{alignat}{3}
    \quad
    & \text{ Dynamical equation } 
    & \
    & \text{ Dispersion }
    \nonumber
    \\
    \nonumber
    \\
    \text{ Classical elasticity: } \quad
    & \rho \partial_t^2 u = N_1 \partial_x^2 u
    & \
    & \rho \omega^2(k) = N_1 k^2
    \\
    \nonumber
    \\
    \text{ Strain gradient: } \quad
    & \rho \partial_t^2 u = N_1\partial_x^2 u - N_2 \partial_x^4 u + \ldots \qquad
    & \
    & \rho \omega^2(k) = N_1 k^2 + N_2 k^4 + \ldots
    \\
    \nonumber
    \\
    \text{ Microinertia: } \quad
    \begin{split}
        & \rho \partial_t^2 u - \rho D_1 \partial_x^2 \partial_t^2{u} +\rho D_2 \partial_x^4 \partial_t^2{u} + \ldots
        \\
        & = N_1 \partial_x^2 u - N_2\partial_x^4 u + \ldots
    \end{split}
    \quad
    & & \rho \omega^2(k) = \frac{N_1 k^2 +N_2 k^4 + \ldots}{1+D_1k^2 +D_2 k^4+\ldots}
\end{alignat}
The dispersion relations have similar behavior at $k^2\to 0$.
The important distinction is the behavior at $k^2\to\infty$.
In this limit, classical elasticity and strain gradient models have the undesirable feature that $\omega^2\to\pm\infty$.
In contrast, microinertia models have a finite limiting value of $\omega$ if we have equal orders of spatial derivatives on both sides of the momentum balance.

We notice that strain gradient models lead to $\omega^2(k)$ that are polynomials, and microinertia models lead to $\omega^2(k)$ that are rational functions.
Rational function approximations have the important advantage that they are more flexible, e.g. they can approximate curves with finite limits at $k^2\to\infty$, while also being able to easily invert to real space and time to give dynamical PDE posed in terms of standard differential operators.

Given this advantage of rational function approximations, our strategy begins by approximating a given exact dispersion relation such as \eqref{disp_rel_periodic_laminate_1d} using a rational function approximation.
Our approximation has the form:
\begin{equation}
\label{general_RF_for_DR}
    \omega_{approx}^2(k^2) = \frac{\sum_{i=1}^N N_i k^{2i}}{D_0 + \sum_{i=1}^N D_i k^{2i}}
\end{equation}
Some important features of this approximation are:
\begin{enumerate}

    \item All coefficients in the rational function approximation are taken to be non-negative. This avoids singularities and the frequencies are real for all $k$.
    
    \item The exponent of the highest-order terms in the numerator and denominator are balanced to ensure a finite limit at $k\to\infty$.
    
    \item The dispersion relation is an even function of $k$, and hence we will use only even powers of $k$ in the rational function approximation. 
    
    \item While our rational function expression only approximately matches the exact given dispersion relation, we constrain the approximation to ensure that we exactly capture the long-wavelength and static behavior. Specifically, we require that our approximation satisfies 
    \begin{equation}
        \lim_{k\to0} \frac{\omega_{approx}^2}{k^2} = c^2
    \end{equation}
    where $c$ is the long-wavelength wave velocity.
    This is ensured by taking $D_0$ as the mean density, and $N_1$ as the effective modulus obtained from static homogenization which is just the harmonic mean of the moduli.
    
    \item We will denote the different rational functions by $RF_{NN}$, where $N$ refers to the exponent of the highest-order term in the numerator or denominator.

\end{enumerate}

\subsection{Homogenized Dynamical Equation and Boundary Continuity Conditions}
\label{sec:1d-1band-jump}

We consider \eqref{general_RF_for_DR} with $N=6$ for concreteness:
\begin{equation}\label{1d_RF66_form}
    \omega_{approx}^2(k^2) = \frac{N_0 {k^2} +N_1 k^4 + N_2k^6}{D_0+ D_1 {k^2} +D_2 k^4+D_3k^6}
\end{equation}
Given an exact dispersion relation obtained from dynamic homogenization of a fine-scale metamaterial model, we select the coefficients $N_0, N_1, N_2, D_0, D_1, D_2$ that best approximate the given exact relation.

We rewrite \eqref{1d_RF66_form} as the dynamical equation in Fourier space:
\begin{equation*}
    \Big(\omega^2(D_0+ D_1 {k^2} +D_2 k^4+D_3k^6)  - (N_0 {k^2} +N_1 k^4 + N_2k^6)\Big) \Tilde{u}(k,\omega) = 0
\end{equation*}
where $\Tilde{u}(k,\omega)$ is the Fourier transform of $u(x,t)$  in both space and time.
We then use the inverse Fourier transform to obtain the (approximate) homogenized dynamical equation in real space and time:
\begin{equation}\label{RF66_inverted_pde_1d_time}
   N_2 \partial_x^6 u  - N_1 \partial_x^4 u  + N_0 \partial^2_x u 
   + D_3 \partial_x^6 \partial_t^2 u 
   - D_2 \partial_x^4 \partial_t^2 u  
   + D_1 \partial_x^2 \partial_t^2 u 
   - D_0 \partial_t^2{u} 
    = 0     
\end{equation}
We emphasize that this Fourier inversion assumes implicitly that the coefficients are all constant in space.

We next construct the corresponding Lagrangian:
\begin{equation}\label{RF66_Lagrangian}
    \begin{split}
        & \calL_{RF_{66}} [u] (t) 
        :=
        \\
        & \int_{x \in \Omega} \Big(N_2 \abs{\partial_x^3 u}^2  + N_1\abs{\partial_x^2u}^2  +
        N_0\abs{\partial_x u}^2 - 
        D_3\abs{\partial_x^3 \partial_t{u}}^2 -
        D_2\abs{\partial_x^2 \partial_t {u}}^2  -  D_1\abs{\partial_x \partial_t{u}}^2 - 
        D_0 \abs{\partial_t{u}}^2 \Big)  \dm x 
    \end{split}
\end{equation}
and the action:
\begin{equation}
\label{eqn:RF66_Action}
    \calS_{RF_{66}} [u] := \int_0^t \calL_{RF_{66}} \dm t
\end{equation}
where $\Omega$ is the body and $[0,t]$ is the time interval of interest.
The variation of $\calS_{RF_{66}}$ with respect to $u(x,t)$ recovers \eqref{RF66_inverted_pde_1d_time} as the Euler-Lagrange equation. 

Here, we make the central assumption of our method.
The dynamical equation \eqref{RF66_inverted_pde_1d_time} was homogeneous in space, i.e. the coefficients $N_0, N_1, N_2, D_0, D_1, D_2, D_3$ were constant.
We now take \eqref{RF66_Lagrangian} and \eqref{eqn:RF66_Action} as the {\em starting points} of our homogenized model, and assume that they hold even if the coefficients are functions of $x$.
Taking the variation of $\calS_{RF_{66}}$ with respect to $u$, we obtain the dynamical equation:
\begin{equation}\label{RF66_Euler_Lag}
\begin{split}
    &-\partial_x^3(N_2(x) \partial_x^3 u)  + \partial_x^2(N_1(x)\partial_x^2u)  - \partial_x (N_0(x)\partial_x u) 
    \\
    & \quad - \partial_x^3(D_3(x)\partial_x^3 \partial_t^2{u}) +\partial_x^2(D_2(x)\partial_x^2\partial_t^2{u}) -  \partial_x(D_1(x)\partial_x \partial_t^2{u})  +D_0(x) \partial_t^2{u} = 0 
\end{split}
\end{equation}
and the corresponding continuity conditions at discontinuities:
\begin{align}
      &\left\llbracket 
        N_2(x) \partial_x^3 u + D_3(x) \partial_x^3 \partial_t^2{u}
     \right\rrbracket  
     =0 \label{RF66_jump_condn_1}
     \\
     &\left\llbracket\partial_x^2 u\right\rrbracket  =0
     \\
      &\left\llbracket
            - \partial_x\left(N_2(x)\partial_x^3 u + D_3(x)\partial_x^3 \partial_t^2{u}\right)
            + \left(N_1(x)\partial_x^2u +D_2(x)\partial_x^2\partial_t^2{u}\right)
        \right\rrbracket =0
     \\
    &\left\llbracket\partial_x u\right\rrbracket  =0
     \\
    &\left\llbracket
        \partial_x^2\left(N_2(x)\partial_x^3 u + D_3(x)\partial_x^3 \partial_t^2{u}\right)
        - \partial_x\left(N_1(x)\partial_x^2u +D_2(x)\partial_x^2\partial_t^2{u} \right)
        + \left(N_0(x)\partial_x u + D_1(x) \partial_x\partial_t^2{u}\right)
    \right\rrbracket 
    = 0
    \\ 
    &\left\llbracket u\right\rrbracket  =0 \label{RF66_jump_condn_6}
\end{align}
The macroscopic homogenized model consists of \eqref{RF66_Euler_Lag} and \eqref{RF66_jump_condn_1}-\eqref{RF66_jump_condn_6}.

We can identify a conserved energy for the homogenized model by multiplying \eqref{RF66_inverted_pde_1d_time} by $\partial_t{u}$ and integrating over $\Omega$:
\begin{equation}\label{RF66_energy_in_time}
\begin{split}
    & E_{RF_{66}}[u](t) :=
    \\ 
    & \frac{1}{2}\int_{x \in \Omega} \left( N_2\abs{\partial_x^3 u}^2  + N_1\abs{\partial_x^2u}^2  + N_0\abs{\partial_x u}^2 + D_3\abs{\partial_t(\partial_x^3  u )}^2+D_2\abs{\partial_t(\partial_x^2u)}^2
    + D_1\abs{\partial_t(\partial_x u)}^2 + 
     D_0 \abs{\partial_t u}^2 \right)  \dm x 
\end{split}
\end{equation}
The energy is clearly positive-definite and its time derivative can be readily computed to be $0$ if there is no energy input from the boundaries of $\Omega$.

\subsection{Numerical Comparisons between the Homogenized and Fine-scale Models}
\label{sec:1d-1band-examples}

We consider the metamaterial with the numerical values of the properties as chosen in Section \ref{sec:num-properties} (Table \ref{tab:composite_prop1d}), and the dispersion relation in \eqref{disp_rel_periodic_laminate_1d}.
The exact dispersion relation and several rational function approximations are shown in Figure \ref{fig:f_vs_K1_1d}, and the values of the coefficients of the rational function approximations are given in Table \ref{tab:RF_coeffs_values}.
The rational function approximations are computed using a least-squares fit.

For both the single- and double-boundary problem, we consider both cases of lamina (a) and lamina (b) being adjacent to the homogeneous material.
While the solution to the corresponding fine-scale models differs slightly, the homogenized model cannot -- even in principle -- capture these minor differences because we do not account for the details of the interface structure.

\begin{figure}[ht!]
    \centering
    \includegraphics[width=0.8\textwidth]{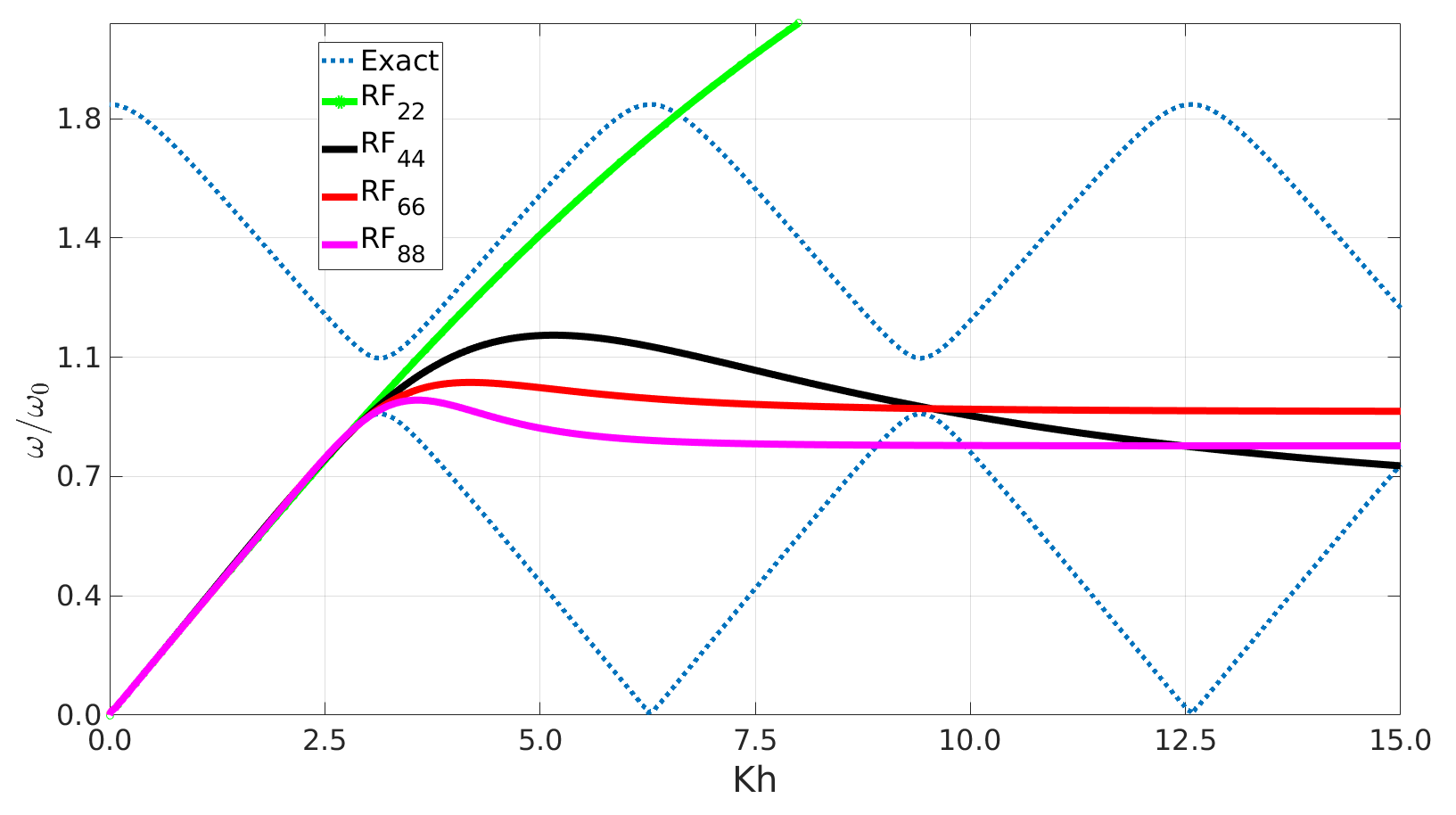}
    \caption{The exact dispersion relation and various rational function approximations.}
    \label{fig:RF_2468_1d_extended}
\end{figure}

\begin{table}[h!]
    \centering
    \begin{tabular}{|c|c|c|c|c|}
        \hline
        \ & $RF_{22}$ & $RF_{44}$&$RF_{66}$ &$RF_{88}$ \\
        \hline
        \hline
        $N_2$ & 0 &\num{1.22451e-02}&  0 & 0 
        \\
        \hline
        $N_3$ & 0 & 0& \num{1.350e-02}&0
        \\
        \hline
        $N_4$ & 0 & 0 & 0 & \num{1.3982e-03}
        \\
        \hline
        $D_1$ & \num{1.100814e-2} & 0 &\num{1.9503e-02}  &\num{2.74576e-04}   \\
        \hline
        $D_2$ & 0 &\num{4.65093e-03}&0&0
        \\
        \hline
        $D_3$ & 0 &0& \num{2.08892e-03}&0
        \\
        \hline
        $D_4$ & 0 & 0&0& \num{2.74576e-04}\\
      \hline 
    \end{tabular}
    \caption{Coefficients in the rational function approximations obtained using  least squares fits to the exact dispersion relation \eqref{disp_rel_periodic_laminate_1d}. Figure \ref{fig:RF_2468_1d_extended} shows plots of the rational functions. $N_1=\num{1.19047619}$ is chosen to match the classical linear elastic limit as $k\to 0$ for all of the approximations. For the homogeneous medium,  $N_1=E_h$ and $D_0=\rho_h$, and all other coefficients are zero. }
    \label{tab:RF_coeffs_values}
\end{table}

\paragraph{Single-Boundary Problem.}

We consider the single-boundary problem described in Section \ref{sec:formulation-problems} (Figure \ref{fig:single_interface_setup}).
The numerical values of the homogeneous material properties are given in Section \ref{sec:num-properties}.
We solve the homogenized model -- \eqref{RF66_Euler_Lag} and \eqref{RF66_jump_condn_1}-\eqref{RF66_jump_condn_6} -- by noticing that \eqref{RF66_Euler_Lag} has piecewise constant coefficients and hence the solution is simply a piecewise superposition of exponentials.
The continuity conditions are then used to obtain the coefficients in the superposition solution.
The solution method is easy and completely standard; the difference is that we need to solve for more coefficients than in the classical case, and we also have more continuity conditions to use.

\begin{figure}[ht!]
    \centering
    \subfloat[Normalized transmitted energy for the fine-scale model and $RF_{66}$ approximation.
    \label{fig:RF_66_energies_1d_1_interface}
    ]
    {\includegraphics[width=0.47\textwidth]{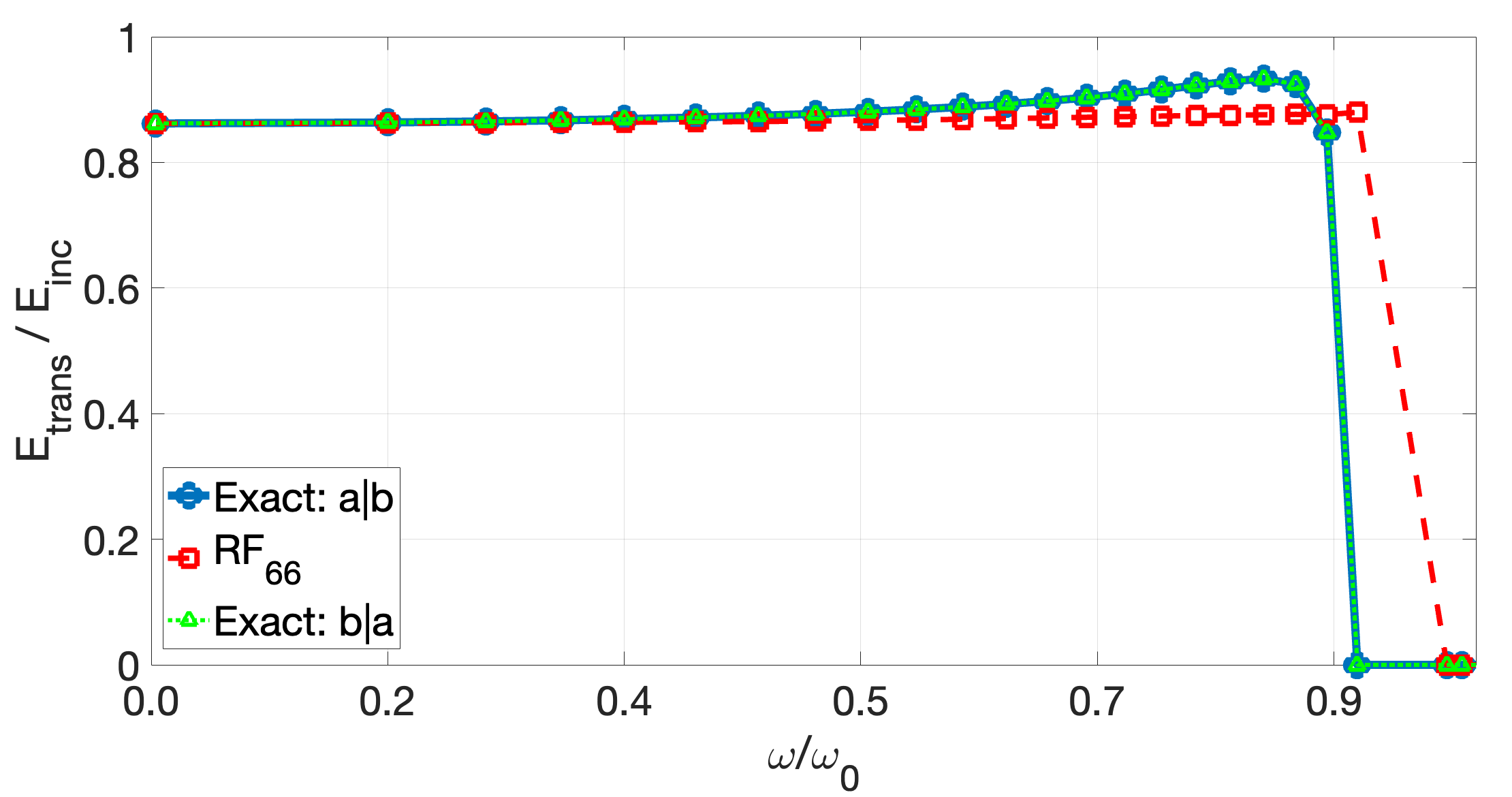}}
    \hfill
    \subfloat[Normalized error in the transmitted energy between the fine-scale model and $RF_{66}$ approximation.
    \label{fig:RF_66_energies_error_1d_1_interface}
    ]
    {\includegraphics[width=0.47\textwidth]{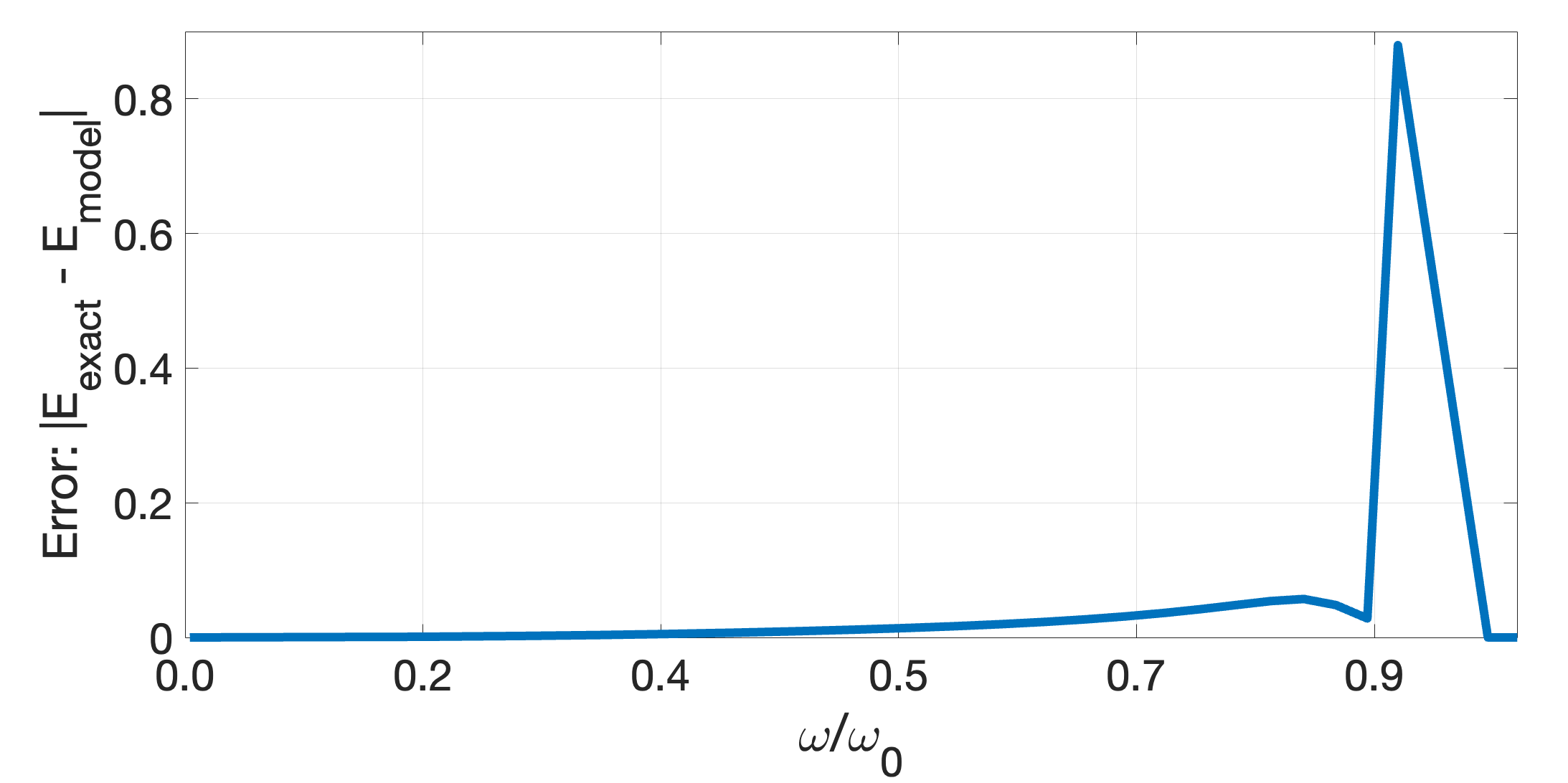}}
    \caption{Transmitted energy predictions and error as a function of frequency for the single-boundary problem in 1-d described in Section \ref{sec:formulation-problems} and Figure \ref{fig:single_interface_setup}. 
    The curve \textit{Exact: a|b} has lamina $(a)$ adjoining the boundary and the curve \textit{Exact: b|a} has lamina $(b)$ adjoining the boundary.}
    \label{fig:RF_66_energies_1d}
\end{figure}

The comparison between the predictions of the exact fine-scale model and the homogenized model is shown in Figure \ref{fig:RF_66_energies_1d}.
Specifically, we show the predictions of normalized energy transmission across the boundary as a function of frequency.
The error in the homogenized model is remarkably small given the significant simplifying assumptions; the error is about $5\%$ at about $80\%$ of the bandgap frequency.

\paragraph{Double-Boundary Problem.}

We next consider the double-boundary problem described in Section \ref{sec:formulation-problems} (Figure \ref{fig:2interface_setup}).
The numerical values of the homogeneous material properties are given in Section \ref{sec:num-properties}.
The solution method for the homogenized problem closely follows the single-boundary case described above.

\begin{figure}[htb!]
    \centering
    \subfloat[Normalized transmitted energy for the fine-scale model and the $RF_{66}$ approximation.
    \label{fig:RF_22_energies_1d_2interface}
    ]
    {\includegraphics[width=0.47\textwidth]{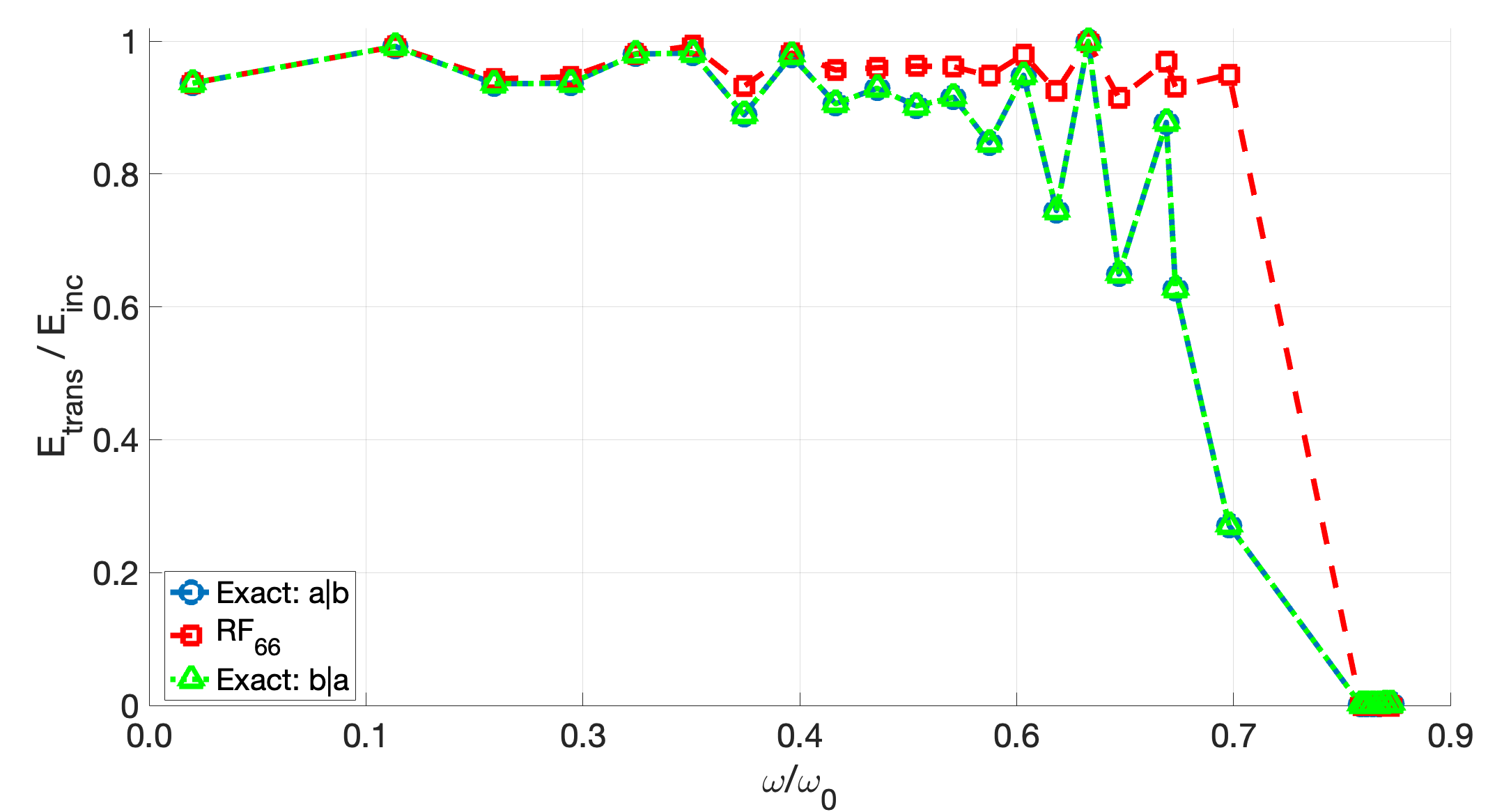}}
    \hfill
    \subfloat[Normalized error in the transmitted energy between the fine-scale model and the $RF_{66}$ approximation.
    \label{fig:RF_22_energies_error_1d_2interface}
    ]
    {\includegraphics[width=0.47\textwidth]{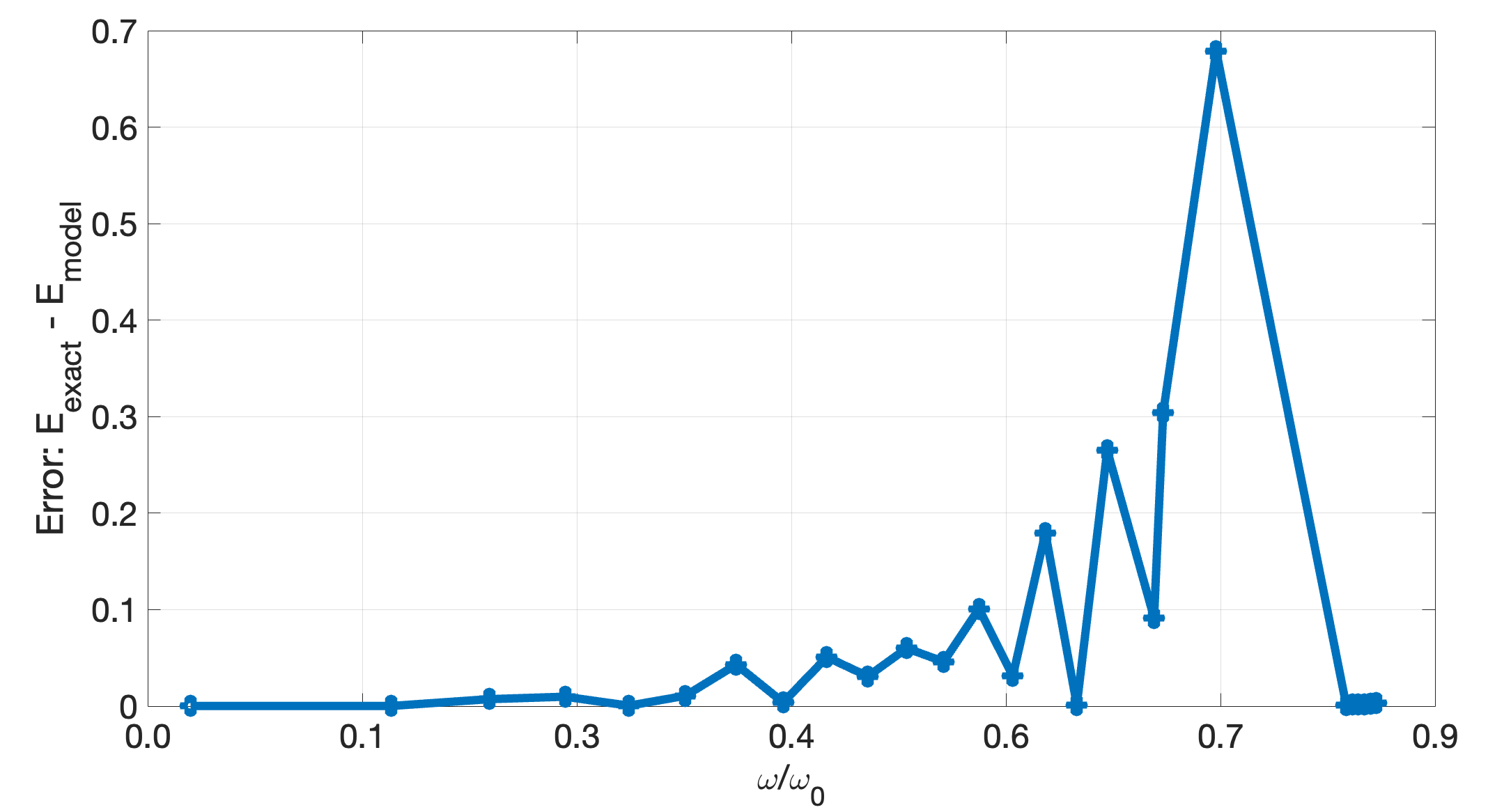}}
    \caption{Transmitted energy predictions and error as a function of frequency for the double-boundary problem in 1-d described in Section \ref{sec:formulation-problems} and Figure \ref{fig:2interface_setup}. We highlight that these predictions are based on the $RF_{66}$ rational function approximation that has a significant overestimate of the bandgap frequency (Figure \ref{fig:RF_2468_1d_extended}), yet it is able to provide good qualitative and fair quantitative predictions.
    }
    \label{fig:RF_22_energies_1d}
\end{figure}

\begin{figure}[htb!]
    \subfloat[$\RE(u)$]{\includegraphics[width=0.5\textwidth]{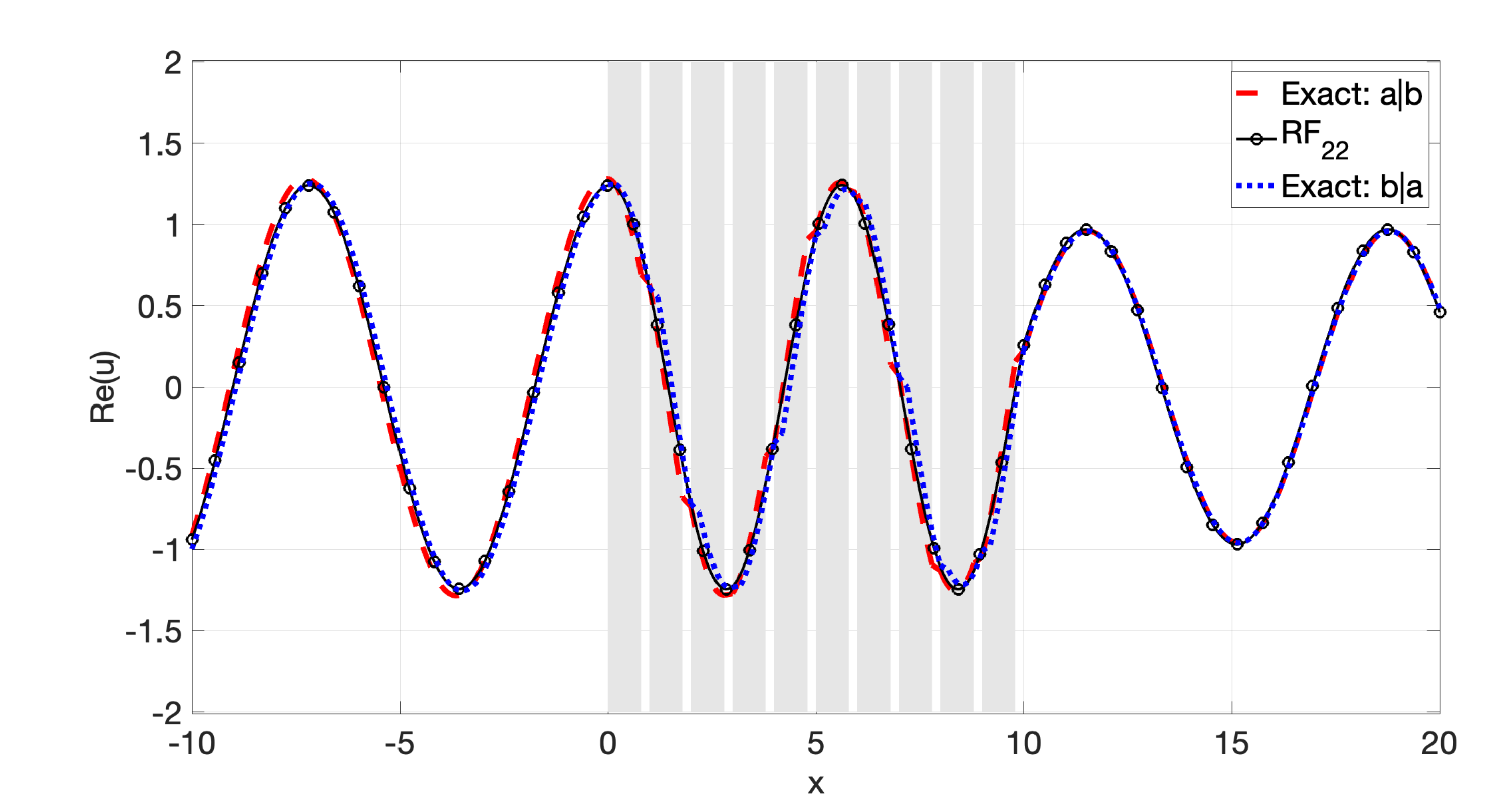}}
	\hfill
	\subfloat[$\IM(u)$]{\includegraphics[width=0.5\textwidth]{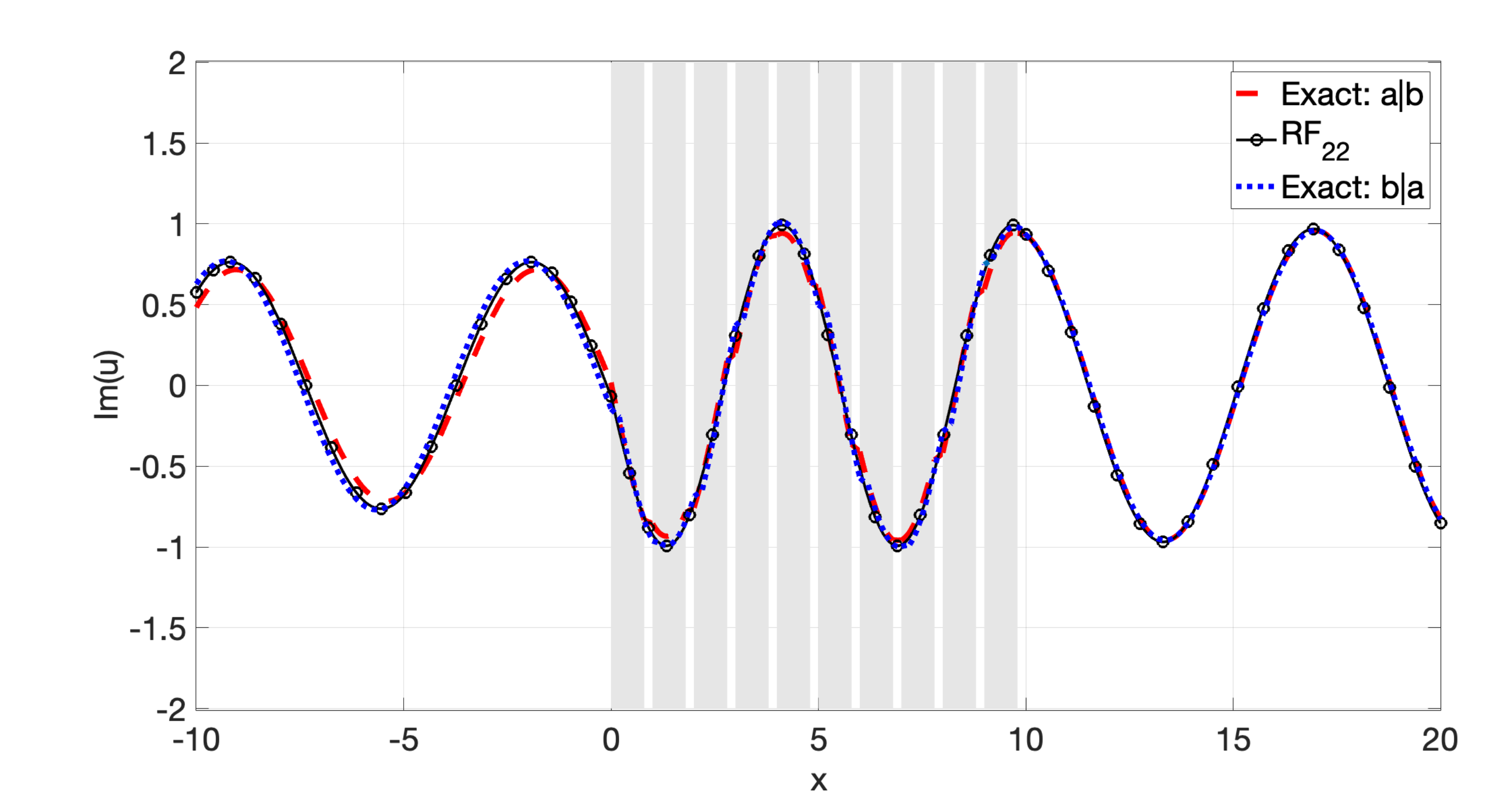}}
	\caption{Comparison between the fine-scale model and the homogenized ($RF_{22}$) model for $\omega/\omega_0 = 0.3552$. The gray shading represents the metamaterial. The curve \textit{Exact: a|b} has lamina (a) adjoining the interface at $x=0$, and similarly curve \textit{Exact: b|a} has lamina (b) adjoining the interface at $x=0$.
	}
	\label{fig:RF_22_f2}
\end{figure}

\begin{figure}[htb!]
	\includegraphics[width=0.8\textwidth]{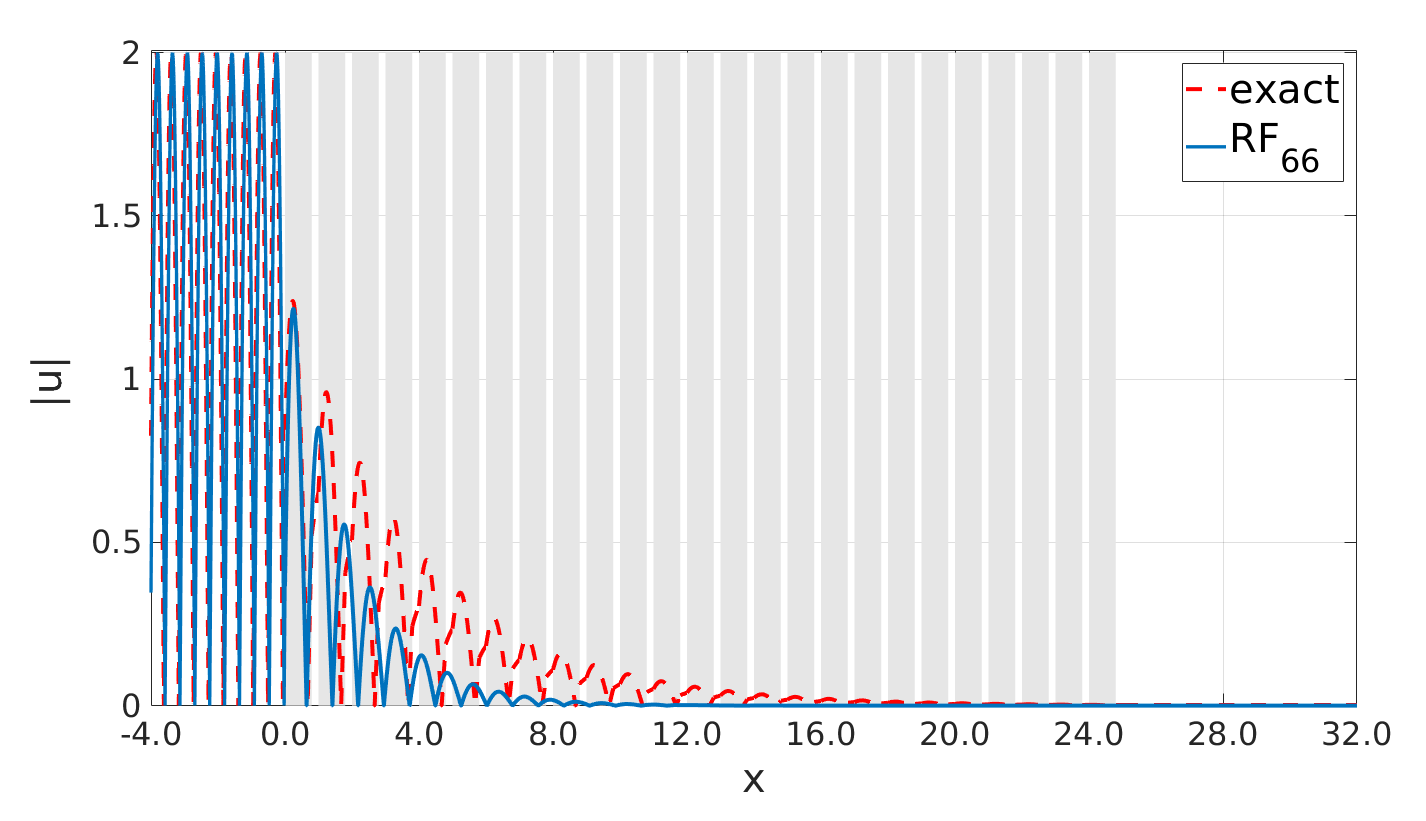}
	\caption{Comparison of the displacement field (absolute value) between the exact fine-scale model and the homogenized $RF_{66}$ approximate model. The frequency $\omega/ \omega_0=1.025$ lies in the bandgap of the metamaterial. The gray shading represents the metamaterial. The homogenized model captures the evanescent waves at the boundary of the metamaterial.}
\label{fig:RF_66_bandgap}
\end{figure}

The comparison between the predictions of the exact fine-scale model and the homogenized model and continuity conditions are shown in Figures \ref{fig:RF_22_energies_1d}, \ref{fig:RF_22_f2} and \ref{fig:RF_66_bandgap}.

Figure \ref{fig:RF_22_energies_1d} shows the transmitted energy as a function of frequency.
We notice the classical dips in the transmitted energy that correspond to the destructive interference between the transmitted waves and the multiple reflected waves.
The homogenized model is able to capture the frequencies at which these dips occur and the overall qualitative structure of the curve well using the $RF_{66}$ approximation which significantly overestimates the frequency near the bandgap (Figure \ref{fig:RF_2468_1d_extended}), with an error of about $20\%$ at about $80\%$ of the bandgap frequency.

Figure \ref{fig:RF_22_f2} shows the displacement field at a frequency about halfway to the bandgap, using the simplest $RF_{22}$ approximation and comparing with the exact fine-scale model.
We find a very good agreement for the real and imaginary parts (i.e., the homogenized model captures the phase information correctly).

If we use a higher-order approximation, e.g. $RF_{66}$, we are able to capture the bandgap much more accurately (Figure \ref{fig:RF_2468_1d_extended}).
Figure \ref{fig:RF_66_bandgap} shows an example of a wave propagating with a frequency $\omega/ \omega_0=1.025$, i.e., it lies in the bandgap of the metamaterial.
We notice that the evanescent modes at the boundary are captured by the homogenized model, though the decay is somewhat faster than that of the exact fine-scale model.

%
\section{Single Band Homogenized Model in Two Dimensions}\label{sec:2d}

We next consider the 2-d setting, and follow the overall strategy of Section \ref{sec:1d}.
For numerical comparisons, we consider the problem described in Figure \ref{fig:perp_interface_setup}.

We begin by approximating the first band of the exact dispersion relation in Figure \ref{fig:3d_RF22_approx} using a rational function approximation of the type:
\begin{equation}\label{RF2d_type2}
    \omega_{approx}^2 = \frac{ \bfN^{(0)}: \bfK\otimes\bfK  + \bfK\otimes\bfK:\bfN^{(1)}: \bfK\otimes\bfK + \dots}{D^{(0)}+\bfD^{(1)} :\bfK\otimes\bfK + \bfK\otimes\bfK:\bfD^{(2)}:\bfK\otimes\bfK + \dots}
\end{equation}
where $\bfN^{(0)}, \bfD^{(1)}$ are second-order tensors, $\bfN^{(1)}, \bfD^{(2)}$ are fourth-order tensors, and so on.
All of these tensors are symmetric and positive definite to avoid singularities in the approximation.

\begin{figure}[htb!]
    \centering
    \includegraphics[width=0.8\textwidth]{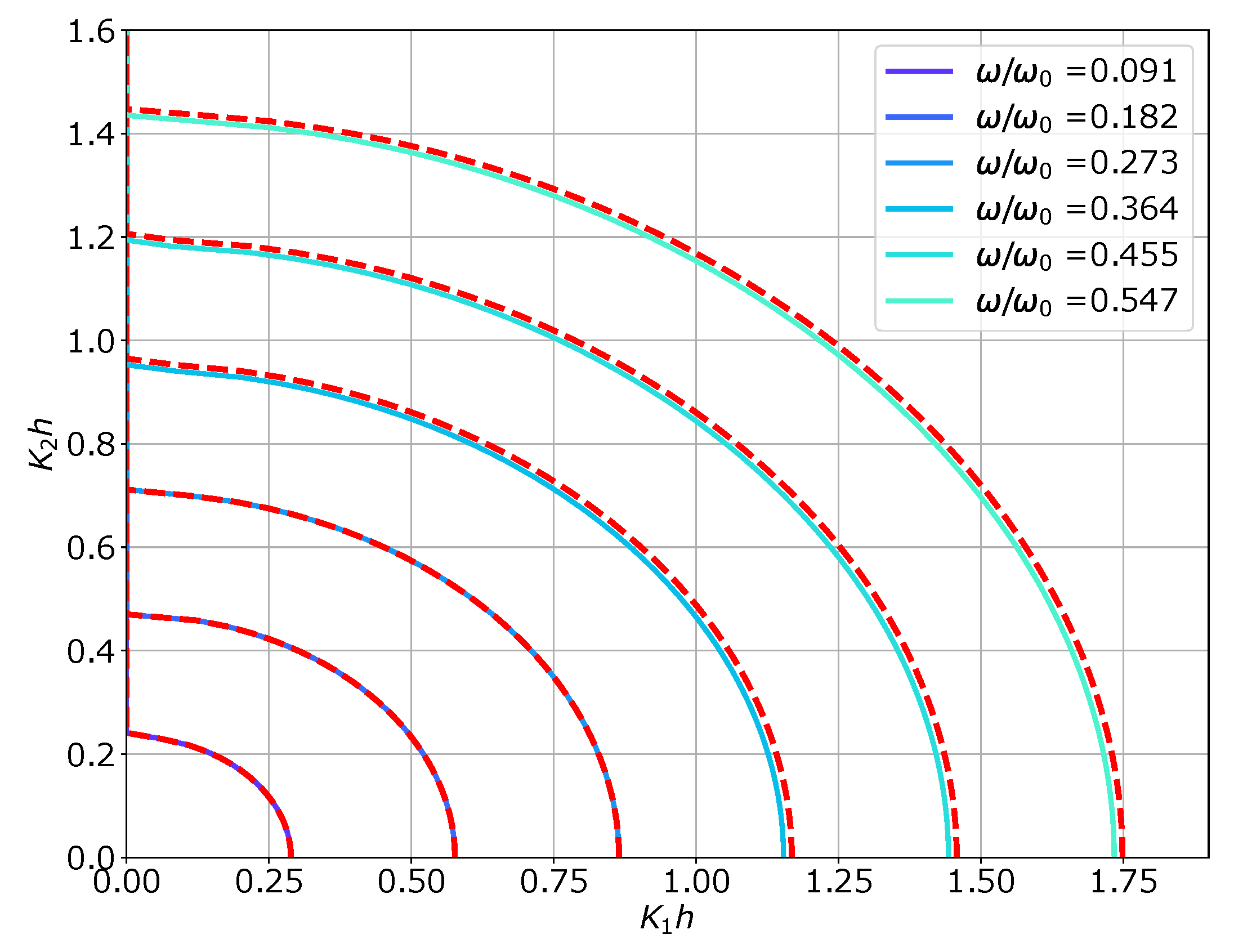}
    \caption{Contours of the exact dispersion relation (solid lines) and its rational function approximation (dashed red lines). }
    \label{fig:K2_K1_f_slice_comparison}
\end{figure}

For the numerical calculations presented in this section, we keep only terms up to $|\bfK|^2$.
A least-squares fit against the exact dispersion relation gives:
\begin{equation}
\label{values_2dRF22}
    D^{(0)} = 1, 
    \quad
    \bfN^{(0)}
    =
    \begin{bmatrix}
        \num{1.344742} & \num{0}
        \\
        \text{sym.} & \num{1.81045806}
    \end{bmatrix},
    \quad
    \bfD^{(1)}
    =
    \begin{bmatrix}
        \num{0.013519456} & 0
        \\
        \text{sym.} & \num{0.01303317}
    \end{bmatrix}
\end{equation}
Figure \ref{fig:K2_K1_f_slice_comparison} compares the exact and the approximate dispersion relations.

\subsection{Homogenized Model and Continuity Conditions}

Transforming \eqref{RF2d_type2} to real space and time gives the 2-d homogenized dynamical equation:
\begin{equation}\label{RF_2d_homogenized_pde}
    - D^{(0)} \partial_t^2 u +  \bfD^{(1)}_{ij}\partial_i \partial_j \partial_t^2 u + \bfN^{(0)}_{ij} \partial_i \partial_j u =0 
\end{equation}

The corresponding Lagrangian $\calL_{2d}$ is:
\begin{equation}\label{RF22_lagrangian_2d}
\begin{split}
    \calL_{2d} = \frac{1}{2}\int_{\Omega} \left(D^{(0)} \abs{\partial_t u}^2 
    - \bfN^{(0)}_{ij}\partial_j u \partial_i u +  \bfD^{(1)}_{ij} \partial_j (\partial_t u) \partial_i (\partial_t u)  \right)\dm \Omega
\end{split}
\end{equation}
Dropping the assumption of homogeneity of the coefficients and using the principle of least action, we get the following continuity conditions on the boundary,
\begin{align}
    \left\llbracket 
        \bfN^{(0)}_{ij}\partial_j u  + \bfD^{(1)}_{ij} \partial_j \partial_t^2 u
    \right\rrbracket \hat{\bfn}_i(\bfx)  &= 0 \label{RF_2d_jump1}
    \\
   \left\llbracket u \right\rrbracket &=0 \label{RF_2d_jump2}
\end{align}
where $\hat{\bfn}(\bfx)$ is the unit normal to the boundary.

\subsection{Numerical Comparisons between the Homogenized and Fine-scale Models}\label{sec:Numerics_2d_1band}

We use the 2-d problem described in Section \ref{sec:formulation-problems} and Figure \ref{fig:perp_interface_setup}, and with the numerical values for the material properties from Section \ref{sec:num-properties} and Table \ref{tab:composite_prop1d}.

The method of solution of the homogenized equations is analogous to the 1-d case.
We first notice that \eqref{RF_2d_homogenized_pde} is linear with piecewise constant coefficients, and hence the solutions are simply  a piecewise superposition of exponential functions.
The coefficients in the superposition are obtained by applying the continuity conditions \eqref{RF_2d_jump1}-\eqref{RF_2d_jump2}.
This gives the scattered displacement field and scattering coefficients for the homogenized model. 

\begin{figure}[ht!]
    \centering
    \subfloat[$\omega /\omega_0=0.29$]{\includegraphics[width=0.5\textwidth]{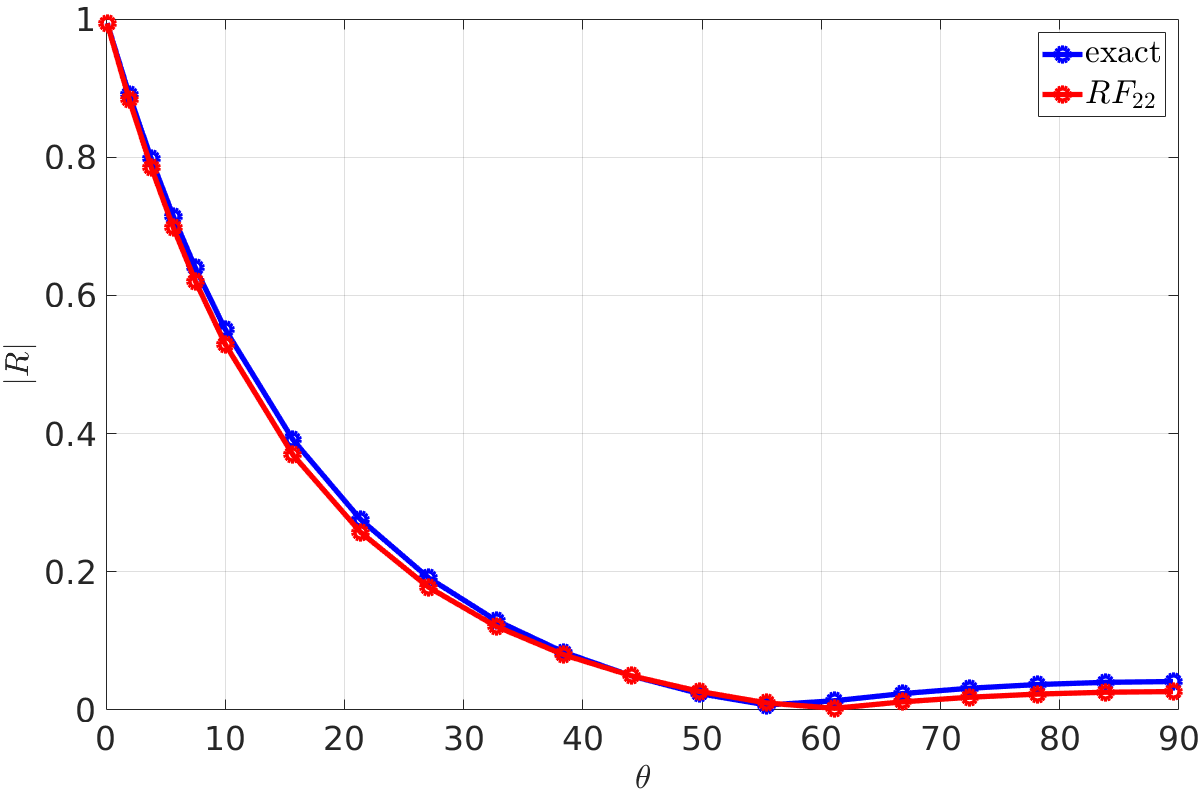}}
    \subfloat[$\omega/\omega_0=0.41$]{\includegraphics[width=0.5\textwidth]{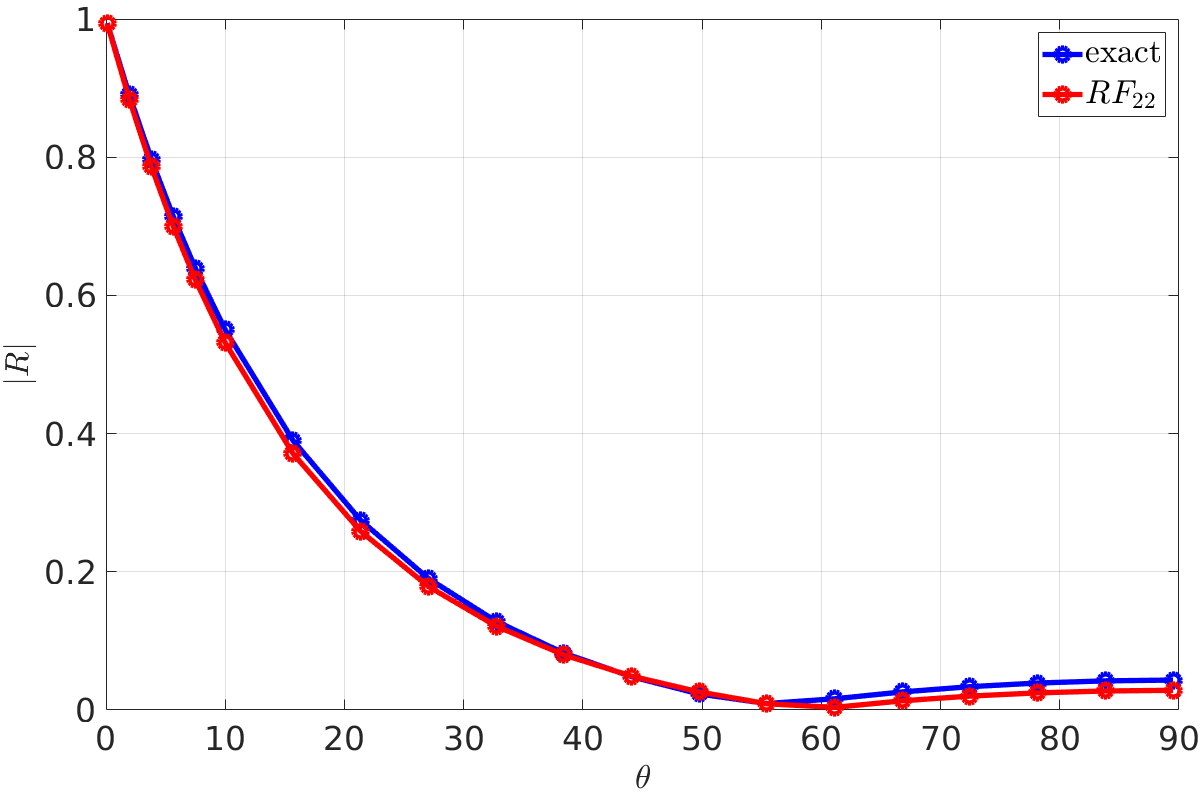}}\\
    \subfloat[$\omega/\omega_0=0.6484$]{\includegraphics[width=0.5\textwidth]{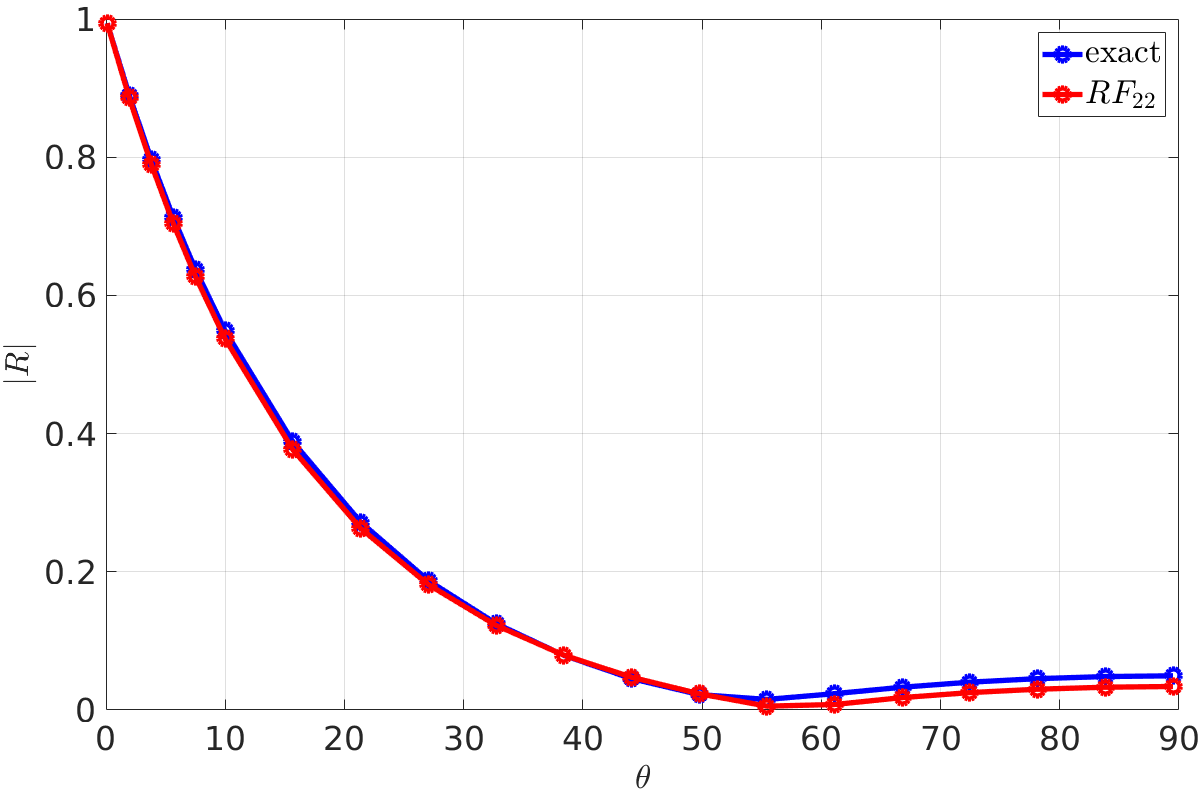}}
    \subfloat[$\omega/\omega_0=1.004$]{\includegraphics[width=0.5\textwidth]{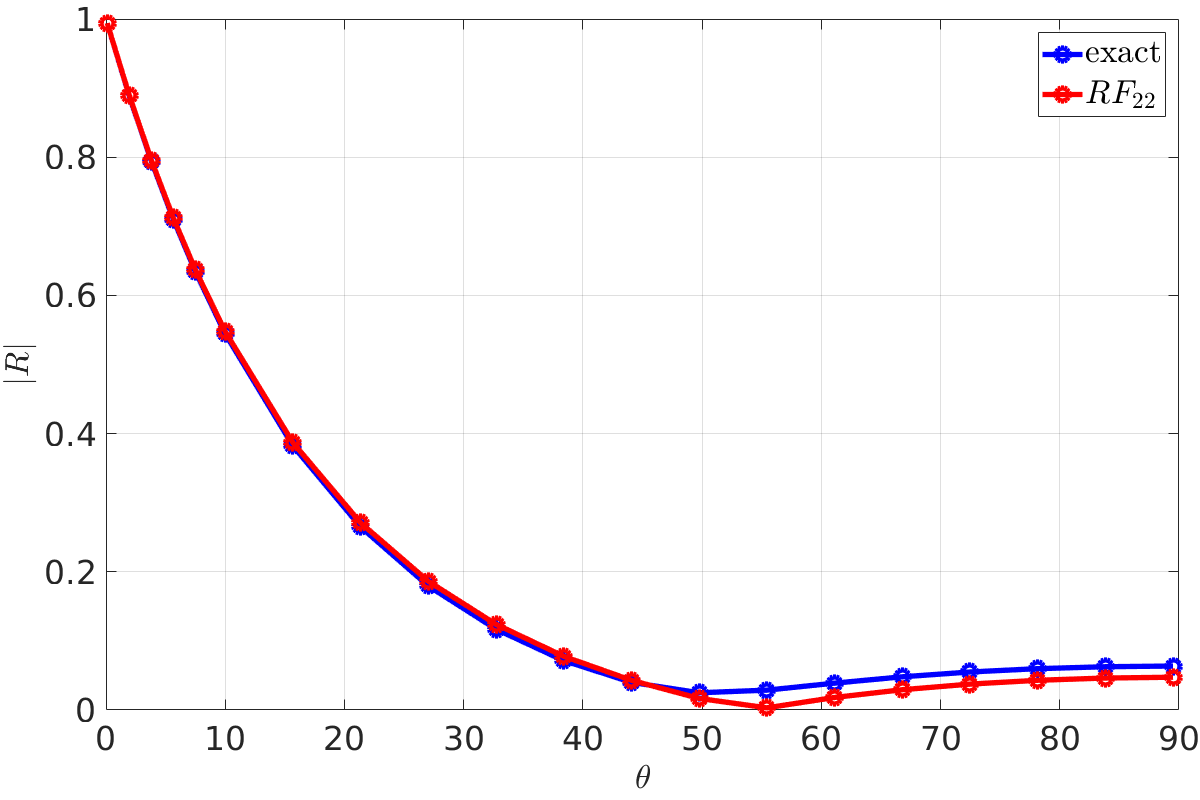}}
    \caption{Comparison between the reflection coefficients predicted by the fine-scale and homogenized models as a function of incident angle for several frequencies.}
    \label{fig:Ref_coeff_comparison}
\end{figure}

Figure \ref{fig:Ref_coeff_comparison} compares the reflection coefficients predicted by the exact fine-scale model and the homogenized model as a function of frequency, for several frequencies and over the entire range of incident angles from grazing to normal incidence.
There is very good quantitative agreement over all the frequencies and incident angles.

\begin{figure}[ht!]
    \centering
    \subfloat[\label{fig:RF22_Re_u_2d}]
    {\includegraphics[width=0.47\textwidth]{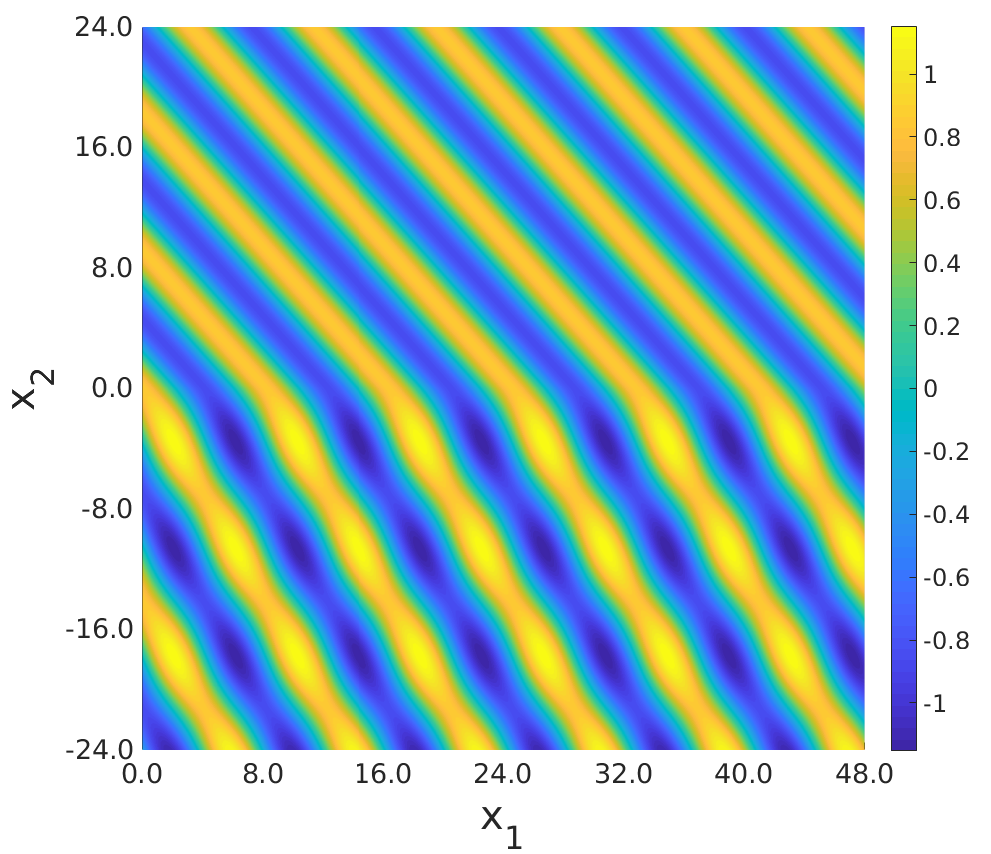}}
    \hfill
    \subfloat[\label{fig:Exact_Re_u_2d}]
    {\includegraphics[width=0.47\textwidth]{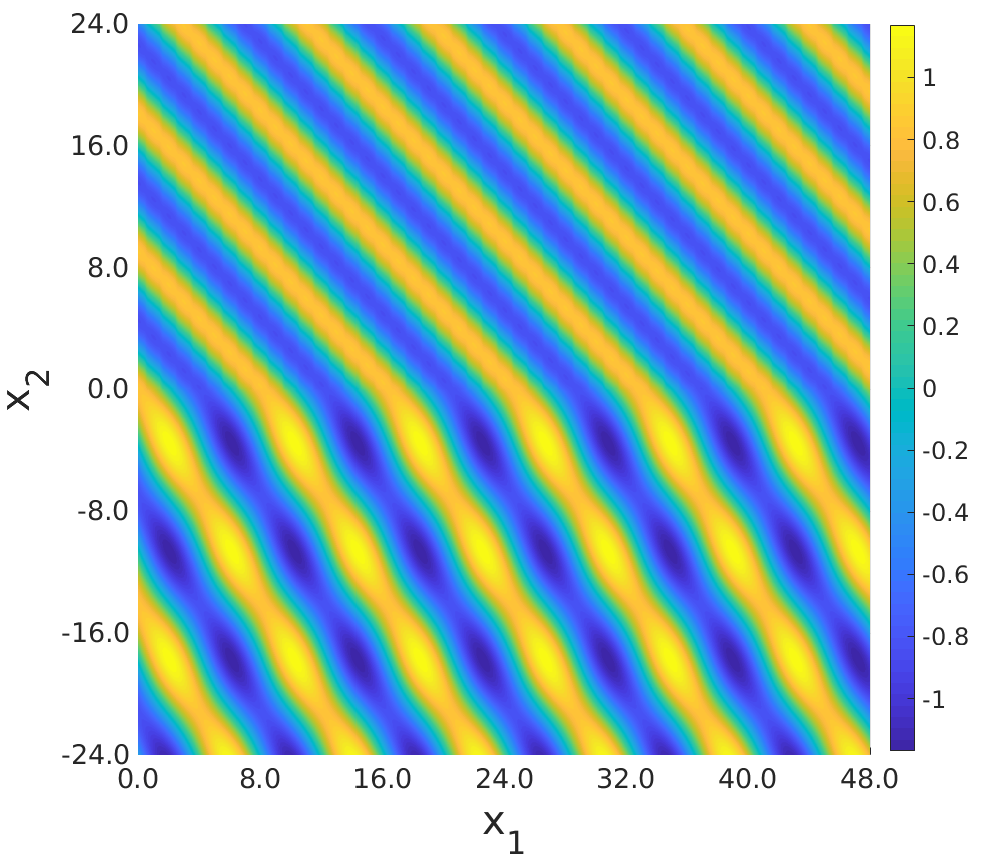}}
    \caption{Real part of the displacement fields in 2-d from the $RF_{22}$ homogenized model (a) and the fine scale model (b). The boundary is along $x_2=0$, the wave is incident at $60^\circ$ to the normal, and $\omega/\omega_0=0.3624$. Visually, the agreement is very good.}
    \label{fig:2d_disp_plot}
\end{figure}

Figure \ref{fig:2d_disp_plot} plots the displacement field of the fine-scale model and the homogenized model for a specific frequency and $\theta=30^\circ$.
Visually, the agreement is very good.

Our consideration of this problem is motivated by \cite{Srivastava2017} that examined the same problem in a different way.
They highlighted that when the metamaterial lamination is normal to the metamaterial boundary, there are a finite number of propagating waves and an infinite number of evanescent waves generated in the metamaterial. 
Standard dynamic homogenization can account only for the propagating waves because it uses the setting of an infinite periodic medium.
If such homogenized models are applied to this problem, the evanescent modes cannot be accounted for, and the classical continuity conditions -- continuity of displacement and traction -- must be satisfied with propagating waves only; this leads to a violation of the energy flux conservation. 
In the case of multiple propagating waves in the $2$-direction, the classical continuity conditions are not sufficient to determine the transmission and reflection coefficients.  

The proposed approach of using higher-order homogenized models gives additional continuity conditions at the metamaterial boundary. 
Further, the higher-order homogenized model supports evanescent waves.
These features address the issues raised in \cite{Srivastava2017} and enable good quantitative predictions.

\section{A Single Homogenized Equation to Describe Two Bands with a Bandgap}
\label{sec:OA}

We now consider the question of modeling multiple bands, and focus on two bands with a bandgap to illustrate the idea in a simple setting.
Prior work, e.g., \cite{Craster2010}, has derived approximations for higher bands, but these approaches are restricted to a single band in the vicinity of an \textit{a priori} selected frequency.
In contrast to such approaches, we obtain a single homogenized equation that is seamlessly valid across the entire frequency range spanning both bands, without having to specify \textit{a priori} the band or limit the frequency range of interest.
Such a model is essential to enable the modeling of complex time-dependent transient loadings, such as shock loading, and waves that are composed of a broad range of frequencies.

Consider two bands, and denote their rational function approximations by $\hat{\omega}_a (k)$ and $\hat{\omega}_b (k)$.
The simplest representation has the form:
\begin{align}
    \omega^2 & =  \frac{nk^2}{1+dk^2} =: \left(\hat{\omega}_a (k)\right)^2 \label{RF22_acoustic_branch}
    \\
    \omega^2 & = \omega_b^2 - \frac{pk^2}{1+qk^2} =: \left( \hat{\omega}_b (k) \right)^2 \label{RF22_optic_branch}
\end{align}
where $n,d,p,q,\omega_b$ are constants.

To obtain a composite dispersion relation that can represent both bands together, we simply use the product of the individual dispersion relations:
\begin{equation}\label{combined_2band}
    (\omega^2 - \hat\omega_a(k)^2)(\omega^2 - \hat\omega_b(k)^2) =0 
\end{equation}
The resulting composite dispersion relation is easy to invert to real space and time to obtain the multiband homogenized equation:
\begin{equation}\label{OA_NLT_inverted_pde}
    \begin{split}
         &\partial_t^4 u - (q+d)\partial_t^4\partial_x^2 u + dq \partial_t^4 \partial_x^4u - 
         (-\omega_b^2dq +pd -nq)\partial_t^2\partial_x^4 u + 
         ( p -\omega_b^2q-\omega_b^2d-n)\partial_t^2\partial_x^2 u
         \\
         & \quad + \omega_b^2\partial_t^2 u 
        -\omega_b^2n\partial_x^2 u 
        + (\omega_b^2nq -np) \partial_x^4 u
                =
                0
    \end{split}
\end{equation}
To simplify notation, we define: $A_1:=q+d, A_2:=qd, A_3:=-\omega_b^2dq +pd -nq, A_4:= p -\omega_b^2q-\omega_b^2d-n, A_5 := \omega_b^2, A_6 := \omega_b^2n, A_7:= \omega_b^2 nq - np$. 

To find the Lagrangian $\calL_{MB}[u](t)$, we first multiply \eqref{OA_NLT_inverted_pde} by a test function $v(x,t)$ and integrate over the spatial domain $\Omega$ and an arbitrary time interval $[t_1, t_2]$.
Using standard integration-by-parts operations, we can simplify to write:
\begin{equation}
\begin{split}
    \int_{x\in\Omega} \int_{t_1}^{t_2} \Big( 
        & \partial_t^2u\cdot\partial_t^2v 
        +  A_1 \partial_t^2\partial_x u \cdot \partial_x \partial_t^2 v 
        + A_2 \partial_t^2\partial_x^2 u \cdot \partial_x^2\partial_t^2 v 
        + A_3 \partial_t\partial_x^2 u \cdot \partial_x^2\partial_t v 
    \\
        & \quad + A_4 \partial_t \partial_x u \cdot \partial_x\partial_t v
        - A_5\partial_t u \cdot \partial_t v
        + A_6\partial_x u \cdot \partial_x v
        + A_7\partial_x^2u \cdot \partial_x^2 v 
    \Big)
    \dm x\dm t
    =
    0
\end{split}
\end{equation}
where we have ignored boundary terms in space and time, by assuming that the initial conditions are given and hence $v$ and its time derivatives are zero at $t=t_1$.

This can be used to write the action:
\begin{equation}
    \calS_{MB} = \int_{t_1}^{t_2} \calL_{MB} \dm t
\end{equation}
where the Lagrangian is:
\begin{equation}
\label{OA_Lagrangian}
    \begin{split}
    &\calL_{MB}
    := \\
    &\int_{x\in \Omega} \Big(
        \abs{\partial_t^2u}^2 
        + A_1\abs{\partial_t^2\partial_x u}^2 
        + A_2\abs{\partial_t^2\partial_x^2 u}^2
        + A_3\abs{\partial_t\partial_x^2 u}^2
        +A_4\abs{\partial_t\partial_x u}^2
        -A_5\abs{\partial_t u}^^2
        + A_6\abs{\partial_x u}^2
        + A_7\abs{\partial_x^2 u}^2
    \Big) \dm x
    \end{split}
\end{equation}

Extremizing $\calS_{MB}$ without assuming that the coefficients are constant gives the homogenized equation:
\begin{equation}
\label{OA_Euler_Lag}
    \partial_t^4u  
    - \partial_t^4 \partial_x \left(A_1 \partial_x u\right) 
    + \partial_t^4 \partial_x^2 \left(A_2 \partial_x^2 u \right) 
    - \partial_t^2 \partial_x^2 \left(A_3 \partial_x^2 u \right) 
    + \partial_t^2 \partial_x \left(A_4 \partial_x u \right) 
    + A_5\partial_t^2 u  
    -\partial_x \left(A_6\partial_x u \right)  
    + \partial_x^2 \left( A_7 \partial_x^2 u \right)
    = 
    0
\end{equation}
 and the continuity conditions:
\begin{align}
\label{OA_jump_conditons-1}
    &\left\llbracket 
        A_1\partial_t^4\partial_x u
        -A_4\partial_t^2 \partial_x u 
        +A_5\partial_t^2 u 
        + A_6\partial_x u
        -\partial_t^4\partial_x \left(A_2 \partial_x^2 u \right)
        +\partial_t^2\partial_x \left(A_3 \partial_x^2 u \right)  
        - \partial_x \left( A_7\partial_x^2 u \right) 
    \right\rrbracket 
    =
    0
    \\
    &\left\llbracket u \right\rrbracket
    =
    0
    \\
    &\left\llbracket 
        A_2\partial_t^4 \partial_x^2 u
        -A_3\partial_t^2 \partial_x^2 u
        + A_7\partial_x^2 u 
    \right\rrbracket
    =
    0
    \\
    &\left\llbracket \partial_x u \right\rrbracket 
    =
    0
    \label{OA_jump_conditons-2}
\end{align}

In Appendix \ref{appendix:energyflux_NLT}, we use $\calL_{MB}$ and $\calS_{MB}$  to obtain the conserved energy \eqref{E_OA_app} and consequently an expression for the energy flux across a surface \eqref{OA_flux_app}.
The energy flux expression is required to compute the transmitted and reflected energies.

\begin{remark}
    An important feature of our homogenized model is the presence of higher-order derivatives in time.
    For a generic higher-order-in-time Lagrangian, Ostrogradsky's theorem \cite{ostrogradski1850malmoires} shows that the corresponding Hamiltonian (energy) has linear instabilities, referred to as {\em Ostrogradsky's instability}. 
    However, our homogenized model is marginally stable in the sense of Lyapunov when the coefficients are chosen in the appropriate range.
\end{remark}

\subsection{Numerical Comparisons between the Homogenized and Fine-scale Models}
\label{sec:Numerics_NLT}

    We use the 1-d single-boundary problem described in Section \ref{sec:formulation-problems} and Figure \ref{fig:single_interface_setup}, with the numerical values for properties listed in Section  \ref{sec:num-properties} (Table \ref{tab:composite_prop1d_bettergap}).

\begin{figure}[htb!]
    \centering
    \includegraphics[width=\textwidth]{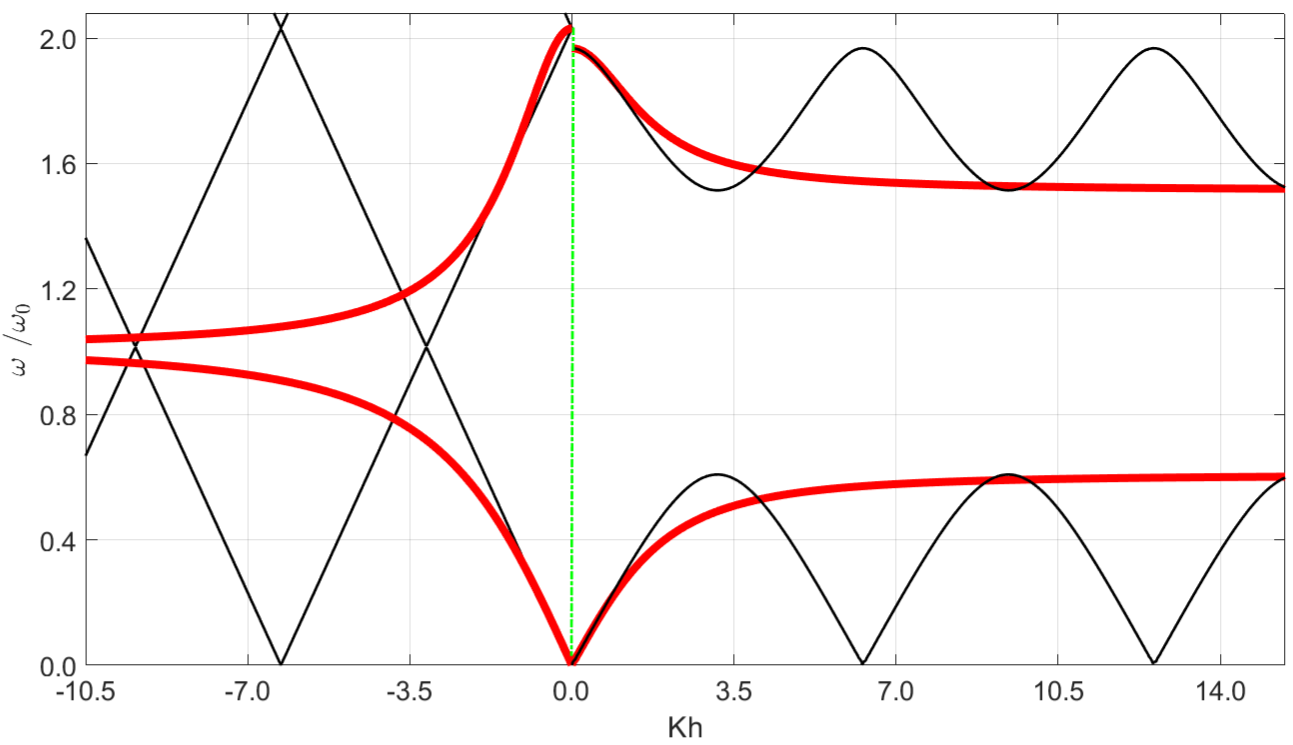}
    \caption{Dispersion relations for the metamaterial ($Kh>0$) and the homogeneous material ($Kh<0$).
    The periodicity of the artificial uniform laminate that represents the homogeneous material is chosen to be the same as the periodicity of the laminate metamaterial.
    The dispersion bands for the metamaterial and the artificially-periodic homogeneous material are shown by the black curves (exact) and red curves (rational function approximation).
    }
    \label{fig:OA_disp_sys_4cases}
\end{figure}

    An important aspect of the homogenized model is that the metamaterial is described using higher-order time derivatives.
    To have a consistent description, we therefore require that the homogeneous material also be modeled using higher-order time derivatives of the same order.
    We achieve this by using considering the homogeneous material as a bilayer laminate metamaterial with identical layers, and constructing a two-band higher-order model (Fig. \ref{fig:OA_disp_sys_4cases}).
    Near $Kh=0$, the linearity of the first band is captured well by the rational function approximation, however, the linearity of the second band is only captured approximately by the second-order rational function approximation.
    Even with this coarse approximation, we find good qualitative and quantitative agreement.

\begin{figure}[htb!]
    \centering
    \includegraphics[width=\textwidth]{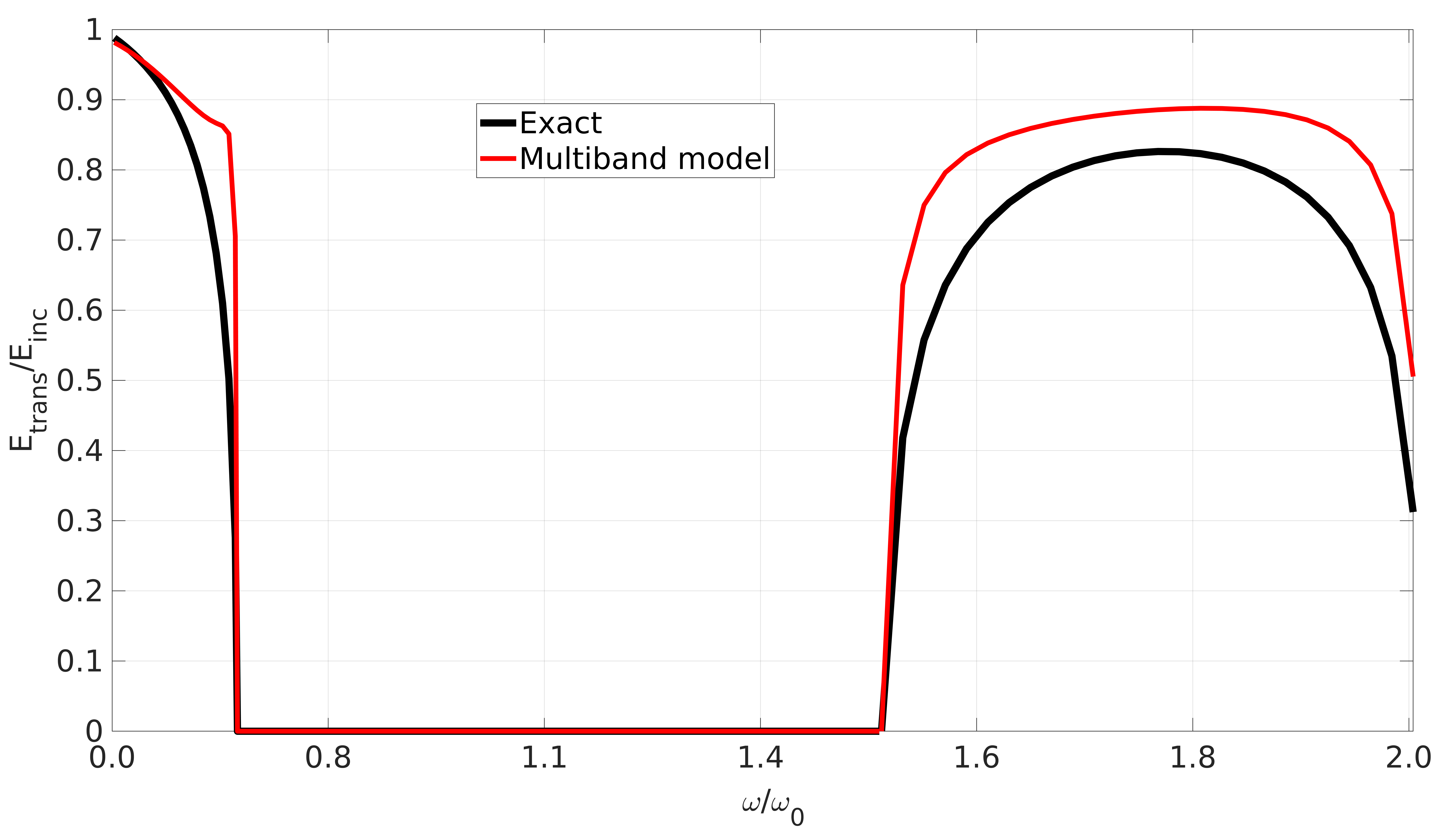}
    \caption{Comparison of the transmitted energies predicted by the exact fine-scale model (black curve) and the homogenized model (red curve) for the approximation shown in Figure \ref{fig:OA_disp_sys_4cases}. The band gap is captured exactly as the rational functions are tuned to match the band gap exactly. }
    \label{fig:OA_energies_4cases}
\end{figure}

    We solve the homogenized model by first using that \eqref{OA_Euler_Lag} is linear with piecewise constant coefficients; hence, the solutions are the standard piecewise superposition of exponentials with unknown coefficients.
    The coefficients are obtained by applying the continuity conditions \eqref{OA_jump_conditons-1}-\eqref{OA_jump_conditons-2}.
    A comparison of the transmitted energy obtained using the fine scale model and the homogenized model is shown in Figure \ref{fig:OA_energies_4cases}. 
    The band gap frequencies can be identified as the range of frequencies over which there is no transmission, and is captured very well by the homogenized model.

\section{A Homogenized Model for Multiple Propagating Waves at High Frequencies}\label{sec:2D-2Band}

\begin{figure}[ht!]
    \centering
    \includegraphics[width=0.8\textwidth]{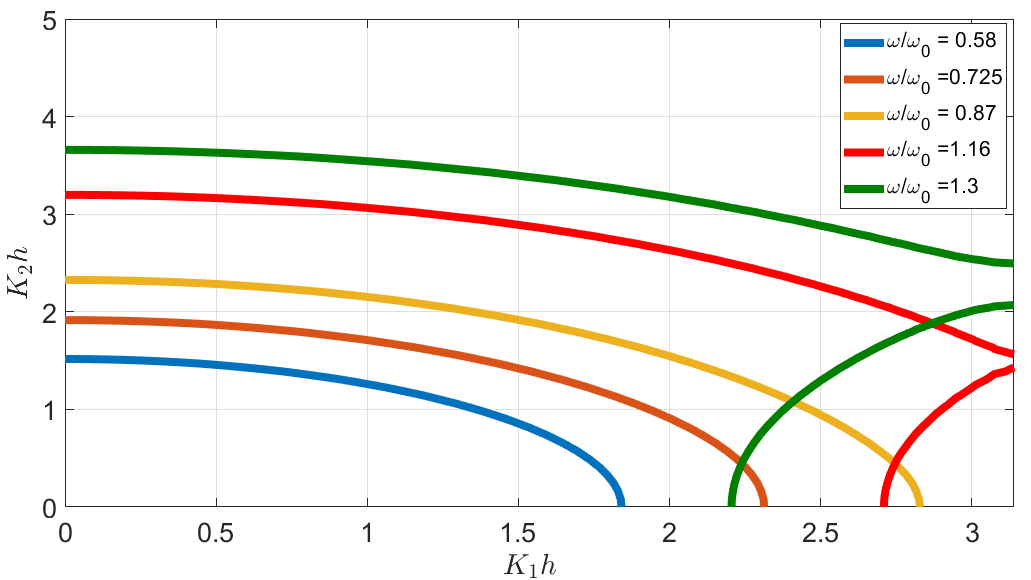}
    \caption{Plot of the dispersion relation as level curves of the frequency as a function of $K_2h$ and $K_1h$. 
    The curves for lower frequencies have only one value of $K_2h$ for every value of $K_1h$, i.e. there is only one propagating wave for that frequency. 
    The curves for frequencies $\omega/\omega_0 = 1.16$ and $\omega/\omega_0 = 1.3$ have two values of $K_2h$ wavevector for some values of $K_1h$, i.e. there are two propagating waves for those frequencies at the corresponding values of $K_1h$.}
    \label{fig:Exact_2d_K1K2}
\end{figure}

As we consider higher frequencies, there are multiple propagating waves in the metamaterial along with infinitely many evanescent waves. 
Figure \ref{fig:Exact_2d_K1K2} shows the $K_1h$ vs. $K_2h$ curves for different frequencies of the first two bands.
We can see that there is a high frequency regime wherein there is a range of values of $K_1h$ for which we have two values of $K_2h$; this implies that there are 2 propagating waves in the metamaterial at those frequencies.

\subsection{Insufficiency of Existing Homogenization Models}

For the problem described in Figure \ref{fig:perp_interface_setup}, where the interface is at right angles to layers in the metamaterial, the displacements fields for the incident, transmitted, and reflected waves are respectively:
\begin{align}
    u_{i} &=
    e^{\imath(\omega t - K_0 \cos{\theta} x_1- K_0 \sin{\theta} x_2)}
    \label{2d_incident_bloch}
    \\
    u_{t} &= \sum_{n=1}^{n=\infty} T_n \Tilde{u}_n (x_1)
    e^{\imath(\omega t - K_1 x_1 - K_2^{(n)}  x_2)}
    \label{2d_ref_bloch}
    \\
    u_{r} &= \sum_{m=-\infty}^{m=\infty} R_m \Tilde{w}_m (x_1)
    e^{\imath(\omega t - K_0 \cos{\theta} x_1 + \kappa_2^{(m)}  x_2)}
    \label{2d_trans_bloch}
\end{align}
This solution uses the Bloch ansatz \cite{Srivastava2017}.
$\Tilde{u}_n(x_1)$, $\Tilde{w}_m(x_1)$ are the Bloch modes in the metamaterial and the homogeneous medium respectively; $K_0$ is the wavevector of the incident wave; and the wavevector in the $x_2$ direction is denoted by $K_2^{(n)}$ for the transmitted field and by $\kappa_2^{(m)}$ for the reflected field.

At low $\omega$, $K_2^{(1)}$ is the only real-valued wavevector, and $K_2^{(n)}$ is imaginary for all $n>1$; this implies that there is only one propagating transmitted wave.

At higher $\omega$ (Fig. \ref{fig:Exact_2d_K1K2}),  $K_2^{(n)}$ for $n=1,2$ are real-valued, and $K_2^{(n)}$ is imaginary for $n>2$; this implies that there are now 2 propagating transmitted waves in the metamaterial.
In general, we could similarly have 1 or more reflected propagating waves in homogeneous medium, but for simplicity we consider only frequencies wherein there is a single reflected propagating wave.

To solve for the scattered propagating waves, we note that $K_1$ is conserved across the interface; i.e., the 2 propagating transmitted waves have the same real value of $K_1h$, but distinct real values of $K_2$.
Even if we completely ignore the evanescent modes, we must determing 3 scattering coefficients: 2 transmission coefficients ($T_1$, $T_2$) for the 2 transmitted propagating  waves, and 1 reflection coefficient ($R_0$) for the reflected propagating wave. 
As discussed in Section \ref{sec:Numerics_2d_1band}, the classical conditions of displacement and traction continuity at the interface provide only 2 equations, but we have to find 3 coefficients.

\subsection{Homogenized Model for Two Bands in Two-dimensions}

We combine the ideas from Section \ref{sec:2d} and Section \ref{sec:OA} to obtain a multiband model for two bands in 2-d. 
The first and second bands of the exact dispersion relation are approximated using rational functions given by $\hat{\omega}_a(\bfK)$ and $\hat{\omega}_b(\bfK)$ respectively.
\begin{align}
    \omega^2 & =  \frac{\bfN : \bfK\otimes\bfK}{1+\bfD:\bfK\otimes\bfK} =: \left(\hat{\omega}_a (\bfK)\right)^2 \label{RF22approx_2d_firstBand}
    \\
    \omega^2 & =  \omega_b^2 - \frac{\bfP : \bfK\otimes\bfK}{1+\bfQ:\bfK\otimes\bfK} =: \left(\hat{\omega}_b (\bfK)\right)^2 \label{RF22approx_2d_secondBand}
\end{align}
where the coefficients $\bfN$, $\bfD$, $\bfP$, and $\bfQ$ are higher-order symmetric tensors. 
The 2-band dispersion relation can be formed by taking the product of the individual approximations, following the idea of \eqref{combined_2band}.

Transforming the dispersion relation to real space and time provides the multiband homogenized dynamical equation in 2-d.
Following prior sections, we construct a Lagrangian corresponding to the dynamical equation, and then use the principle of least action to derive the dynamical equation and the corresponding jump conditions:
\begin{equation}\label{2d_2band_inverted_pde}
    \begin{split}
         &\partial_t^4 u 
         - \partial_i(\bfA^{(1)}_{ij}\partial_t^4\partial_j u) 
         + \partial_{mn}(\bfA^{(2)}_{mnpq}\partial_t^4 \partial_{pq}u) 
         - \partial_{mn}(\bfA^{(3)}_{mnpq}\partial_t^2 \partial_{pq}u) 
         + \partial_i(\bfA^{(4)}_{ij}\partial_t^2\partial_j u) 
         + A^{(5)}\partial_t^2 u
         \\
         & \quad  -  \partial_i(\bfA^{(6)}_{ij}\partial_j u) 
         + \partial_{mn}(\bfA^{(7)}_{mnpq} \partial_{pq}u) 
         =
         0
    \end{split}
\end{equation}
and
\begin{align}
\label{2d_2Band_EffectiveJumps-1}
    &\left\llbracket 
        -\partial_t^4\bfA^{(1)}_{mn}\partial_n u
        +\partial_t^4\partial_{n}(\bfA^{(2)}_{mnpq} \partial_{pq} u) 
        - \partial_t^2\partial_{n}(\bfA^{(3)}_{mnpq} \partial_{pq} u)
        + \partial_t^2(\bfA^{(4)}_{mn} \partial_{n} u)
        + \partial_{n}(\bfA^{(7)}_{mnpq} \partial_{pq} u)
        \right. \nonumber
        \\
        &\left.
        - (\bfA^{(6)}_{mn} \partial_{n} u)\right\rrbracket \hat{\bfn}_{m}
    =
    0
    \\
    &\left\llbracket u \right\rrbracket
    =
    0
    \\
    &\left\llbracket 
        -\bfA^{(2)}_{mnpq}\partial_t^4 \partial_{pq} u
        +\bfA^{(3)}_{mnpq}\partial_t^2 \partial_{pq} u
        - \bfA^{(7)}_{mnpq} \partial_{pq} u
    \right\rrbracket\hat{\bfn}_n
    =
    0
    \\
    &\left\llbracket \partial_i u \right\rrbracket 
    =
    0
    \label{2d_2Band_EffectiveJumps-4}
\end{align}
Here, $\hat{\bfn}(\bfx)$ is the unit normal to the boundary.

For conciseness, we have defined: $\bfA^{(1)}_{ij}:=\bfQ_{ij}+\bfD_{ij},
\bfA^{(2)}_{mnpq}:=\bfD_{mn}\bfQ_{pq},
\bfA^{(3)}_{mnpq}:=-\omega_b^2\bfD_{mn}\bfQ_{pq} +\bfD_{mn}\bfP_{pq} -\bfN_{mn}\bfQ_{pq},
\bfA^{(4)}_{ij}:= -\omega_b^2(\bfQ_{ij}+\bfD_{ij}) + \bfP_{ij}- \bfN_{ij},
\bfA^5 := \omega_b^2,
\bfA^{(6)}_{ij} := \omega_b^2\bfN_{ij},
\bfA^{(7)}_{mnpq}:=\omega_b^2\bfN_{mn}\bfQ_{pq} - \bfN_{mn}\bfP_{pq}$.

\subsection{Quantitative Results} \label{sec:numerics_2d_2band}

\begin{figure}[ht!]
    \centering
    \includegraphics[width=0.8\textwidth]{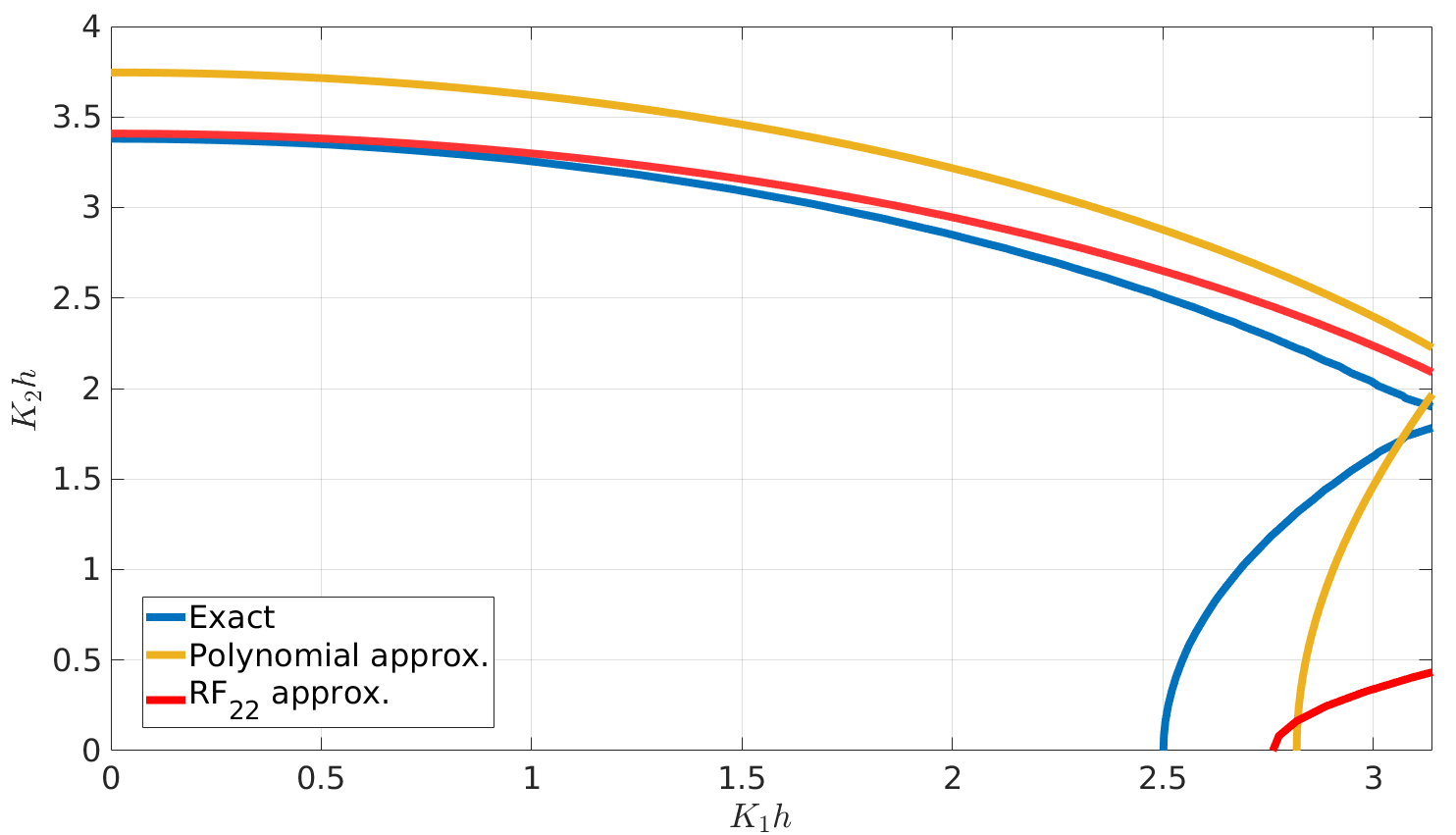}
    \caption{$K_2h$ vs. $K_1h$ plot for $\omega/\omega_0 = 1.218$: Approximation of the exact dispersion relation (solid blue curve) using rational functions (dotted red curve) and polynomials (dotted black curve). }
    \label{fig:2d_2band_K1K2}
\end{figure}

We begin by constructing the homogenized approximation of the metamaterial.
To enable comparisons, we compare the rational function construction with a polynomial construction of the dispersion relation (Fig. \ref{fig:2d_2band_K1K2}).
Using the simple rational functions of the form shown in \eqref{RF22approx_2d_firstBand} and \eqref{RF22approx_2d_secondBand}, the coefficients are obtained from a least-squares fit to the exact dispersion relation:
\begin{equation} \label{2d_2band_RF22_coeffs}
    \begin{split}
        &\bfN^{(0)}
        =
        \begin{bmatrix}
            \num{1.344742} & \num{0}
            \\
            \text{sym.} & \num{1.81045806}
        \end{bmatrix},
        \quad
        \bfD^{(1)}
        =
        \begin{bmatrix}
            \num{0.013519456} & 0
            \\
            \text{sym.} & \num{0.01303317}
        \end{bmatrix},
        \\
        \omega_b = 6.399088, \quad 
        &\bfP^{(0)}
        =
        \begin{bmatrix}
            \num{22.535058} & \num{0}
            \\
            \text{sym.} & \num{0}
        \end{bmatrix},
        \quad
        \bfQ^{(1)}
        =
        \begin{bmatrix}
            \num{0.83561577} & 0
            \\
            \text{sym.} & \num{1.5551709}
        \end{bmatrix}
    \end{split}
\end{equation}

Consider the 2-d problem described in Section \ref{sec:formulation-problems} and Figure \ref{fig:perp_interface_setup}. 
Choosing a frequency $\omega/\omega_0 = 1.218$, at which we see two propagating transmitted waves in the metamaterial, we consider an incident wave with wavevector $K_1h=2.9$.
Since $K_1h$ is conserved across the interface, the metamaterial has two propagating waves with $K_1h=2.9$ but different values for $K_2h$.
The homogeneous medium has just one propagating reflected wave.
The fine-scale energy transmission is computed following Section 5 from \cite{Willis2016}; the transmitted energy ($E_{trans}/E_{inc}$) from the fine-scale model is $0.7914$, evaluated by using as many evanescent modes as needed for the error to be below $0.1\%$.

The multiband homogenized model provides enough continuity conditions \eqref{2d_2Band_EffectiveJumps-1}-\eqref{2d_2Band_EffectiveJumps-4} to solve for the scattering coefficients. 
The transmitted energy computed from the homogenized model based on the polynomial approximation of the exact dispersion relation is $0.2844$, while the model based on the rational functional approximation predicts a transmitted energy of $0.7001$.

We draw 2 conclusions: first, the homogenized model based on rational function approximations does well quantitatively given the relative crudeness of the approximation; and, second, the comparison to the polynomial-based approximation suggests that a good approximation of the dispersion relation provides significant gains in accuracy.
Based on these conclusions, we conjecture that further refinement of the dispersion approximation, through the use of higher-order models and the corresponding higher-order continuity conditions, will provide further quantitative accuracy.

\section{Nonlocal Limit Models}
\label{sec:integral_models}

Section \ref{sec:OA} showed the emergence of 4th-order time derivatives when we considered 2 bands.
In principle, it is straightforward to consider as many bands as one wishes by continuing the approach in \eqref{combined_2band}.
For instance, if we wish to consider $N$ bands, we can write each band in terms of a rational function approximation $\hat\omega_i(k)$, and construct the composite dispersion relation as:
\begin{equation}
\label{eqn:limit-dispersion}
    \prod_{i=1}^N \left(\omega^2 - \hat\omega_i^2(k)\right) = 0
\end{equation}
The resulting dispersion relation will include (even) powers of $\omega$ of up to order $\omega^{2N}$, and the corresponding homogenized model in real time and space will include  (even) time-derivatives of up to order $2N$.

We recall that the fine-scale model will generically have an infinite number of bands.
This raises the question of what limit homogenized model we might expect if we aim to represent all bands in our homogenized model.
We could think of each higher-order time derivative as containing information from further out in time; the (heuristic) limit of an infinite order time-derivative might then correspond to a nonlocal operator in time.

Our heuristic reasoning parallels work in peridynamics \cite{Silling2000}, wherein formal arguments based on Taylor expansions have been used to relate nonlocal integral operators in space to an infinite series of spatial derivatives of increasing order, following the idea of matching dispersion relations \cite{Weckner2005,dayal2017leading,breitzman2018bond,seleson2009peridynamics,seleson2016consistency}.
From that perspective, one could think of peridynamics as a special case of the limit models discussed here, in that it is the limit of a single band model.

In \cite{dayal2017leading}, it was shown formally that the limit of microinertia models -- corresponding to a single-band model in the perspective put forward in this paper -- is consistent with a nonlocal integral equation of the form:
\begin{equation}
\label{eqn:peridyn-KNL}
    \rho \partial_t^2{u}(x,t) 
    + \int_{x^\prime\in\Omega} 
    C_1(x,x^\prime) \left(\partial_t^2{u(x,t)} - \partial_t^2{u(x^{\prime},t)}\right)\dm x^{\prime} 
    =
    \int_{x^\prime \in \Omega} 
    C_2(x,x^\prime) \left({u(x^\prime,t) - {u(x,t)}}\right)\dm x^{\prime}
\end{equation}
where $C_1(x,x^\prime)$ and $C_2(x,x^\prime)$ are spatial kernels. 
We note that the continuity conditions for nonlocal operators are much weaker than in models with higher-order derivatives.

Formally, if we consider the limit of time derivatives of higher order, we can arrive at a limit model that has nonlocality in both space and time is of the form:
\begin{equation}\label{full_non_local_EOM}
\begin{split}
    &\left(u(x,t-T_0)-2u(x,t)+u(x,t+T_0)\right)
    \\
    +
    &\int_{x^\prime\in\Omega } 
    \int_{\tau=0}^{\tau=T_0}
    C_1(x,x^\prime)
    K(\tau)
    \Big( 
      u(x,t-\tau)-2u(x,t)+u(x,t+\tau)
      -u(x^\prime,t-\tau)+2u(x^\prime,t)-u(x^\prime,t+\tau)
    \Big) 
    \dm \tau \dm x^\prime
    \\
    & = 
    \int_{x^{\prime}\in \Omega}
    C_2(x,x^\prime)
    \left(u(x^\prime,t) - u(x,t)\right)
    \dm x^{\prime}
\end{split}
\end{equation}
$K(\tau)$ is a time-kernel, and $T_0$ is the fixed extent of non-locality in time.
The structure of the term that is nonlocal in both space and time follows from \eqref{eqn:peridyn-KNL}, and the term that is nonlocal in time alone has a Taylor expansion that matches the structure of the terms in \eqref{eqn:limit-dispersion}.

\begin{figure}[htb!]
    \centering
    \includegraphics[width=0.5\textwidth]{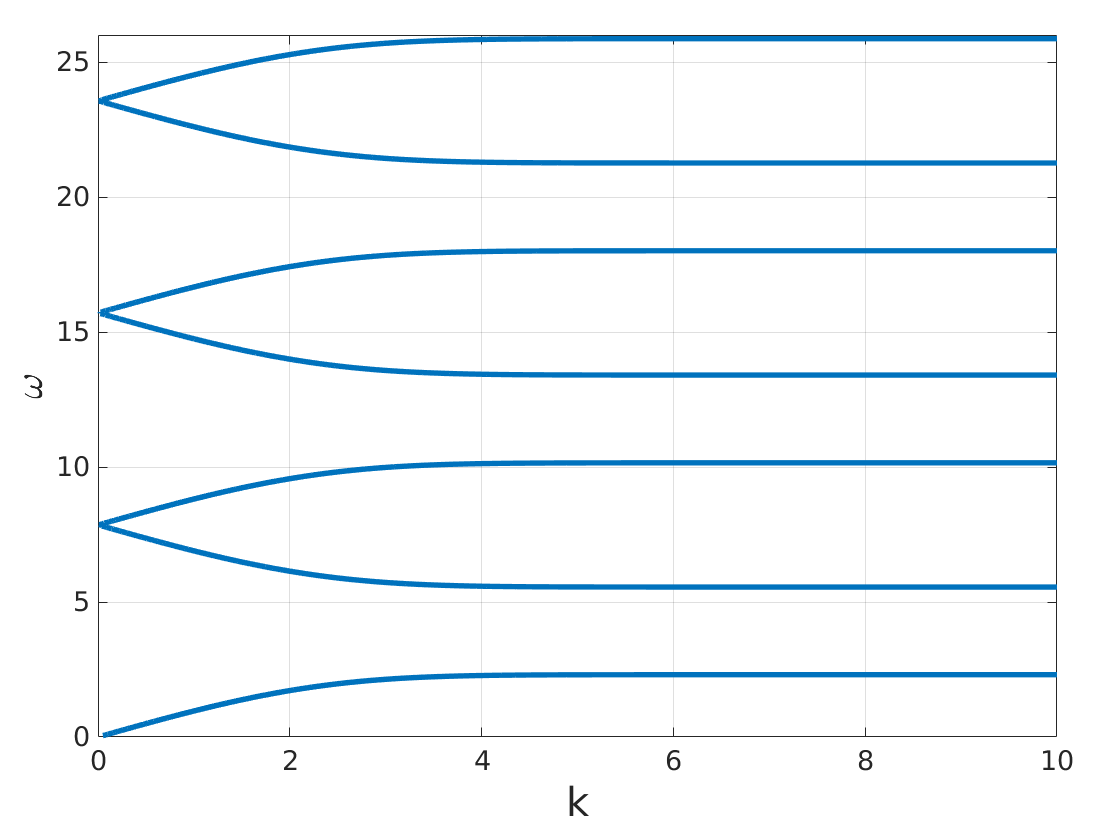}
    \caption{Dispersion relation of the space-time non-local model from \eqref{full_NLT_DR}.}
    \label{fig:full_NLT_DR}
\end{figure}

To examine the bandstructure of \eqref{full_non_local_EOM}, as an example we choose $T_0 = 0.8,  K(\tau) = 1$, $C_1(y) = 0.01\frac{e^{- y^2}}{\sqrt{\pi}}$, $C_2(y)=2 \frac{e^{- y^2}}{\sqrt{\pi}}$ gives the dispersion relation:
\begin{equation}\label{full_NLT_DR}
    \begin{split}
        \frac{\sin{\omega}}{50\omega}
        - \frac{e^{\frac{-k^2}{4}}\sin{\omega}}{50 \omega}
        - \frac{25\sin{(0.4\omega)}^2}{8}
        - \frac{99 e^{\frac{-k^2}{4}}}{50}
        + \frac{99}{50}
        =
        0
        \end{split}
\end{equation}
Figure \ref{fig:full_NLT_DR} shows the dispersion relation: it appears physically reasonable and has an infinite number of dispersion bands.

\begin{remark}[Willis-type Nonlocal Homogenized Descriptions]
    Willis' approach \cite{Willis1980,willis1981variational} provides an alternative strategy to develop a nonlocal limit model; it is based on using Fourier space and averaging over space and time.
    In contrast, while we also work using Fourier space, we then homogenize over space and time through the mechanism of considering the dispersion bands.
\end{remark}

\section{Discussion}

We have presented an approach to develop homogenized models for metamaterials that are posed in space and time, and consequently potentially applicable to general settings such as complex geometries and complex time-dependent loading and to a broad range of frequencies.
The homogenized models have higher-order derivatives in space and time, leading to nonstandard continuity conditions at the boundaries between metamaterials.

While the predictions of the homogenized models compares well with exact fine-scale models, there are a number of important open questions.
The most central question is a justification of the method in a more ``bottom-up'' fashion.
While our approach provides models that are much simpler to use than the rigorous approaches, e.g. \cite{Cornaggia2020}, it is important to understand the connection between our approach and these rigorous approaches; in particular, if it is possible to view our approach as a simplified version of the rigorous approaches under certain assumptions.
Similarly, computational homogenization schemes developed for dynamic situations, e.g. \cite{liu2017variational,hui2014high,pham2013transient,sridhar2018general} provide alternative, but potentially more expensive, approaches to the problems considered here.
These may provide a route to connect the approach in this paper to finer scales using computational approaches; however, it is also important to demonstrate that these methods can effectively capture the evanescent waves at interfaces, which appears to be open at this time.

A second open question is how to perform practical numerical calculations with our model, given the nonclassical higher-order derivatives.
In the context of finite element (FE) methods, two options appear feasible.
First, second-order rational functions provided fairly good accuracy.
They have spatial derivatives of 4th-order, making them amenable to Isogeometric FE methods which provides sufficient continuity for such derivatives \cite{kamensky2019tigar}.
Second, mixed FE methods are a possible approach for problems that require very smooth interpolations, though mixed methods can require care to formulate correctly.
The implementation of numerical methods will enable application of the approach to numerous problems of current interest in the broad area of metamaterials.

An important avenue for future development is to go beyond problems with a scalar unknown field to consider vector problems that come up in elasticity and electromagnetism.
While the overall principles do not change, one must account for the tensorial nature of problems with vector unknowns.
In particular, the dispersion relation that must be inverted now requires tensorial quantities to be approximated as functions of the wavevector; a component-wise approximation is perhaps a feasible strategy to attack this problem.
Similarly, we have restricted our work in this paper to 1 and 2 spatial dimensions and this must be extended to 3 dimensions; however, in our view, this is a relatively straightforward next step.

We highlight that this work, and the many others cited here, are all in the linear regime.
Homogenization of wave scattering in the nonlinear regime is largely unexplored, though we mention a study of wave scattering at an interface between a phase-transforming solid and a linear elastic solid \cite{abeyaratne1992reflection}.
On the other hand, we also note that the linear dispersion behavior is important even in nonlinear problems, for instance to govern the effective dissipation \cite{dayal2006kinetics} or stability \cite{abeyaratne2014macroscopic}.

We note that some aspects of our overall approach -- namely, the approximation of dispersion relations and subsequent inversion to real space and time -- is broadly similar to strategies used in the construction of nonreflecting boundary conditions (NRBC) \cite{bayliss1980radiation, givoli1990non}.
In the NRBC literature, such an approach has been used to construct wave equations that propagate waves in a tailored manner.

In the context of linear homogenization, we note other interesting works that obtain homogenized models that have features in common with our homogenized models.
Specifically, the elimination of fine-scale degrees of freedom has been found to give rise to effective nonlocality in space by \cite{pratapa2018bloch,drugan1996micromechanics}, and also to give rise to very nonstandard dynamical behavior in \cite{milton2007modifications,paola2012exact}.
Further, \cite{Willis1980,willis1981variational,tartar1991memory,antonic1993memory,tartar1989nonlocal} find the emergence of nonlocality in time due to homogenization.
Our work provides heuristic insight into an alternative mechanism: specifically, it suggests the emergence of temporal nonlocality as an outcome of spatial homogenization through the consideration of higher-order dispersion bands that corresponding to modes with finer spatial structure.

A consequence of our approach is also the perspective that peridynamics can be considered as the limit of a spatially averaged model that consider only a single dispersion band; the physical origin of the nonlocality in space is the coarse-graining of the fine structure.
Analogously, a richer model that considers all dispersion bands will provide the full spatio-temporally nonlocal model that is equivalent to averaging in time and space, and can be useful for recent work in nonlocal modeling that includes wave propagation and dynamics \cite{you2020data,you2022data}.



\begin{acknowledgments}
    We thank George Gazonas, Robert Lipton, Graeme Milton, and Luc Tartar for useful discussions; NSF (1635407), ARO (W911NF-17-1-0084, MURI W911NF-19-1-0245), ONR (N00014-18-1-2528), AFOSR (MURI FA9550-18-1-0095), and the Science and Engineering Research Board Fellowship from the Department of Science and Technology of the Indian Government for financial support; and Pittsburgh Supercomputing Center for computing resources.
    This paper draws from the doctoral dissertation of Kshiteej Deshmukh at Carnegie Mellon University.
\end{acknowledgments}

\appendix
\makeatletter
\renewcommand*{\thesection}{\Alph{section}}
\renewcommand*{\thesubsection}{\thesection.\arabic{subsection}}
\renewcommand*{\p@subsection}{}
\renewcommand*{\thesubsubsection}{\thesubsection.\arabic{subsubsection}}
\renewcommand*{\p@subsubsection}{}
\makeatother

\section{Energy Flux for the Nonlocal-in-Time Model}
\label{appendix:energyflux_NLT}

Given a Lagrangian with higher-order time derivatives, an expression for the corresponding Hamiltonian can be obtained by considering the time derivative of the Lagrangian. 
We first show this for a single particle system with an equation of motion given by:
\begin{equation}\label{particle_EOM}
    \partial_t^4 u_p + C \partial_t^2 u_p  - u_p =0, 
\end{equation}
where $u_p \equiv u_p(t)$ describes the particle motion and $C$ is a given constant.  
The corresponding Lagrangian for which \eqref{particle_EOM} is the Euler-Lagrange equation is:
\begin{equation}
\label{particle_Lagrangian}
    \calL_p[u_p](t) =   \frac{1}{2} \left(\abs{\partial_t^2u_p}^2  - C\abs{\partial_t u_p}^2 - \abs{u_p}^2\right) 
\end{equation}
Taking the time derivative of \eqref{particle_Lagrangian}, we get:
\begin{equation}
\begin{split}
    \deriv{\calL_p}{t} 
    =  
    & \partial_t^2 u_p \cdot \partial_t^3 u_p 
    - C \partial_t u_p \cdot \partial_t^2 u_p 
     - u_p \cdot \partial_t u_p
    \\
    =
    & \partial_t u_p \cdot \underbrace{\left( \partial_t^4 u_p + C \partial_t^2 u_p  - u_p \right)}_{=0 \text{ from } \eqref{particle_EOM}}
    + \deriv{\ }{t} \left( \abs{\partial_t^2 u_p}^2  
    - C\abs{\partial_t u_p}^2 - \partial_t u_p \cdot \partial_t^3 u_p \right)
\end{split}
\end{equation}
Therefore, we have the conserved energy:
\begin{equation}
    E_{p}[u_p(t)] := 
    \abs{\partial_t^2 u_p}^2  - C\abs{\partial_t u_p}^2 - \partial_t u_p \cdot \partial_t^3 u_p - \calL_p
    = 
    \frac{1}{2}\abs{\partial_t^2 u_p}^2 -
    \frac{1}{2}C\abs{\partial_t u_p}^2
    - \partial_t u_p \cdot \partial_t^3 u_p + \frac{1}{2} \abs{u_p}^2
\end{equation}
We apply a similar approach to the Lagrangian in \eqref{OA_Lagrangian}.
We take the time derivative of \eqref{OA_Lagrangian}, and use the chain rule, integration-by-parts, and \eqref{OA_Euler_Lag} to obtain terms over the bulk and over the boundary.
The bulk terms give the conserved energy:
\begin{equation}\label{E_OA_app}
    \begin{split}
        &E_{MB}[u(x,t)] := 
        \\
        2\int_{x\in \Omega} \Bigg( & \frac{1}{2}\abs{\partial_t^2{u}}^2 
        - \frac{A_5}{2} \abs{\partial_t{u}}^2
        - \partial_t{u} \cdot \partial_t^3 u
        +\frac{A_1}{2}\abs{\partial_t^2{(\partial_xu)}}^2
        + \frac{A_4}{2} \abs{\partial_t{(\partial_xu)}}^2
        - A_1 \partial_t{(\partial_xu)} \cdot \partial_t^3(\partial_xu) \Bigg.
        \\
        &
        \quad + \Bigg.\frac{A_2}{2}\abs{\partial_t^2{(\partial_x^2u)}}^2 
        + \frac{A_3}{2} \abs{\partial_t{(\partial_x^2u)}}^2
        - A_2  \partial_t{(\partial_x^2u)}\cdot \partial_t^3(\partial_x^2u)  
        - \frac{A_6}{2}\abs{\partial_xu}^2
        - \frac{A_7}{2}\abs{\partial_x^2u}^2  \Bigg)
        \dm x 
    \end{split}
\end{equation}
and the boundary terms give an expression for the flux through a surface:
\begin{equation}\label{OA_flux_app}
\begin{split}
    E_{{MB}_{flux}}[u(x,t)] := 
    &
    2 \left( \partial_t u \cdot A_1\partial_x\partial_t^4u - 
    \partial_t u \cdot  \partial_x(A_2\partial_x^2\partial_t^4u)
    + \partial_x\partial_t u \cdot A_2\partial_x^2\partial_t^4u
    + \partial_t u \cdot  \partial_x(A_3\partial_x^2\partial_t^2u)\right.
    \\
    & \quad
    \left.- \partial_x\partial_t u \cdot  A_3\partial_x^2\partial_t^2u
    - \partial_t u \cdot  A_4\partial_x\partial_t^2u
    + \partial_t u \cdot  A_6\partial_x u
    + \partial_x\partial_t u \cdot  A_7\partial_x^2u
    - \partial_t u \cdot  \partial_x(A_7\partial_x^2u)
    \right)
\end{split}
\end{equation}

\newcommand{\etalchar}[1]{$^{#1}$}

\end{document}